\documentclass[final,1p,times]{elsarticle}

\usepackage{amssymb}
\usepackage{amsthm}
\usepackage{amscd}
\usepackage{amsmath}
\usepackage{amsfonts}
\usepackage{amssymb}
\usepackage{graphicx}
\usepackage{color}
\usepackage{framed}
\usepackage{comment}
\usepackage{chngcntr}
\numberwithin{equation}{section}
\newtheorem{theorem}{Theorem}[section]
\newtheorem{proposition}[theorem]{Proposition}
\newtheorem{lemma}[theorem]{Lemma}

\newtheorem{definition}[theorem]{Definition}
\newtheorem{remark}[theorem]{Remark}
\usepackage{mathrsfs}
\usepackage{titletoc}
\biboptions{sort&compress}
\usepackage[breaklinks,colorlinks]{hyperref}
\usepackage[pagewise,mathlines]{lineno}
\usepackage{geometry}
\usepackage{bm}

\makeatletter

\newcommand{\Rmnum}[1]{\expandafter\@slowromancap\romannumeral #1@}
\makeatother

\usepackage{bbding}

\journal{***}

\begin{document}

\begin{frontmatter}

\title{Fractional Laplace operator on finite graphs}
		
\author[author1]{Mengjie Zhang}
\ead{zhangmengjie@mail.tsinghua.edu.cn}
\author[author2]{Yong Lin}
\ead{yonglin@mail.tsinghua.edu.cn}
\author[author3]{Yunyan Yang{\footnote{corresponding author}}}
\ead{yunyanyang@ruc.edu.cn}		
\address[author1]{Department of Mathematical Sciences, Tsinghua University, Beijing 100084, China}	
\address[author2]{Yau Mathematical Sciences Center, Tsinghua University, Beijing 100084, China}
\address[author3]{School of Mathematics, Renmin University of China, Beijing, 100872, China}	

\begin{abstract}

Nowadays great attention has been focused on the discrete fractional Laplace operator as the natural counterpart of the continuous one. In this paper, we discretize the fractional Laplace operator $(-\Delta)^{s}$ for an arbitrary finite graph and any positive real number $s$. It is shown that $(-\Delta)^{s}$ can be explicitly represented by eigenvalues and eigenfunctions of the Laplace operator  $-\Delta$. Moreover, we study its  important properties, such as
$(-\Delta)^{s}$ converges to   $-\Delta$ as $s$ tends to $1$; while $(-\Delta)^{s}$ converges to the identity map as $s$ tends to $0$ on a specific function space.  For related problems involving the fractional Laplace operator, we consider the fractional Kazdan-Warner equation and obtain several existence results via variational principles and the method of upper and lower solutions.

\end{abstract}

\begin{keyword}
fractional Laplace operator; finite graphs; Kazdan-Warner equation; variational method
\\
\MSC[2020] 35R02; 35A15
\end{keyword}
		
\end{frontmatter}
	
\section{Introduction}

As a core example of a nonlocal pseudo differential operator,
fractional  Laplace operator has been a classical topic in functional and harmonic analysis all along in Euclidean space.
We highly recommend some classical books and papers that discuss this theme in detail, such as
 \cite{A1, A, A78,A12}  and the many references given therein.
The fractional Laplace operator and the corresponding  nonlinear problems are experiencing impressive applications in various fields,
such as water waves   \cite{A32,A42,A35},
optimization   \cite{A37},
crystal dislocation  \cite{A8,A47},
conservation laws   \cite{A9},
singular set of minima of variational functionals   \cite{A69,A55},
stratified materials  \cite{A81},
minimal surfaces   \cite{A20},
gradient potential theory  \cite{A71,T8}
and so on.

For the definition of the fractional Laplace operator on a continuous function space,  readers can refer to  \cite{Kwa15}, which collects ten equivalent definitions. Now we focus on one of these definitions, specifically the definition given by the heat semigroup.
To be exact,  fixing the fractional exponent $0<s<1$,  the  representation of the fractional Laplace-Beltrami  operator  $(-\Delta)^s $  is given by
\begin{eqnarray} \label{ee}
 (-\Delta)^s u(x)=\frac{s}{\Gamma(1-s)} \int_0^{+\infty}
\left(u(x)-e^{t \Delta}u(x)\right) t^{-1-s} d t, ,\ \forall u\in C^{\infty}(\mathbb{R}^N),
\end{eqnarray}
where    $ \Gamma(\cdot)$ is the Gamma function and $e^{t \Delta}$
denotes the heat semigroup of the  Laplacian  $\Delta$.
Noting that the Laplace-Beltrami operator is a geometric terminology denoting the negative Laplace operator, we slightly misuse names and do not distinguish between the Laplace-Beltrami operator and the Laplace operator in this paper.
One of foremost properties for $(-\Delta)^{s}$ is the asymptotic behavior as  ${s\rightarrow 1^-}$ and ${s\rightarrow 0^+}$. In view of (\cite{A}, Proposition 4.4), the following statement holds:
   let $ N > 1$,  there are for any $u\in C_0^\infty({\mathbb{R}^N})$,
\begin{eqnarray}\label{a} \lim_{s\rightarrow 1^-} (-\Delta)^{s} u=-\Delta u\ \text{ and } \ \lim_{s\rightarrow 0^+} (-\Delta)^{s} u=u .\end{eqnarray}
Moreover, note that $ (-\Delta)^s$ is a nonlocal linear operator, but the operator $ -\Delta$ is a local linear operator, which is the essential difference between them. However,   $(-\Delta)^{s}$ and  $-\Delta$   are   closely related as $s$ tends to $1$.

 In the last few years, interest in the discrete fractional Laplace operator has emerged.
 There are some woks about the discrete form of the Euclidean fractional Laplace operator  \eqref{ee}  for any $s\in(0,1)$.
  On the integer lattice  graph  $\mathbb{Z}$,
  \eqref{ee}  is  studied in  \cite{CRSTV2018} and     \cite{KN23},  and  an equivalent form is given as
  \begin{align}\label{Wang}
(-\Delta)^s u(x)
=   \sum_{y \in  \mathbb{Z}} K_s(x-y) \big(u(x)-u(y)\big)  ,\ \forall u\in L^2(\mathbb{Z}),
\end{align}
where the discrete kernel $K_s$ is given by
  \begin{align*}
  K_s(x)= \frac{s}{\Gamma(1-s)} \int_0^{+\infty}e^{t \Delta}1_0(x) t^{-1-s} dt .
\end{align*}
and $1_y$   denotes the characteristic function of $y\in \mathbb{Z}$ with $1_0(x-y)= 1_y(x)$.
On $N$-dimensional  lattice graph  $\mathbb{Z}^N$ ($N\geq 1$),    \eqref{ee}  is  studied in \cite{LR2018}. On a stochastically complete graph $G = (V, E, \mu,  {w})$  with standard measure and weight, \eqref{ee}  is  studied in \cite{Wang-frac}   and an equivalent form is wrote as  \eqref{Wang}. The stochastically complete graph   means that  the weighted graph $G$ satisfies
\begin{align}\label{heat-3} \sum_{y\in V} p(t, x, y)\mu(y)=1, \ \forall t >0,\  \forall x\in V,\end{align}
where $p(t,x,y)=e^{t \Delta}(1_y(x)/\mu(y))$ denotes the heat kernel on $G$.
As for other equivalent forms of the Euclidean fractional Laplace operator, there are also corresponding discrete versions and related problems, see for examples   \cite{CRSTV2015, Tarasov, CRSTV2017, LM2022, JZ22, XZ20}.

 In addition, nonlinear partial differential equations on graphs are also a very attractive topic.
As we know,
the Kazdan-Warner equation arises from the basic geometric problem of prescribing Gaussian curvature on the Riemann surface.
  Grigor'yan-Lin-Yang \cite{A-Y-Y-2} first studied the corresponding finite graphs version
of the Kazdan-Warner equation, and obtained characterizations of the solvability of this equation.
Later, the existence of a solution in the critical case was proved by Ge  \cite{Ge}.
 More recently,  Liu-Yang \cite{L-S} considered this equation with multiple solutions on finite graphs;
 the Kazdan-Warner equation is  studied by  Ge-Jiang \cite{Ge-Jiang}  on infinite graphs, and by Keller-Schwarz \cite{KS18} on canonically compactifiable graphs;
 Sun-Wang \cite{S-w}  studied Kazdan-Warner equations  via degree theory;
 for other related works, we refer the readers to  \cite{A-Y-Y-3,  H-W-Y,    L-Y2,   H-X, H-S, HS, H-L-W}  and references therein. \\

  One aim of this paper is to   provide a discrete version of
the fractional Laplace operator  $(-\Delta)^s$  in \eqref{ee}  for  any fractional  exponent $ s>0$  on an arbitrary finite graph, which is inspired by 
 \eqref{Wang}.
 Another aim of this paper is  to study  a nonlinear problem involving the discrete operator $(-\Delta)^s$,  namely,  the solvability of the fractional Kazhdan-Warner equation.

 To be specific,  let  $G = (V, E, \mu, w)$ be a connected finite graph.
Our discussion is presented in two cases: $0<s<1$ and $s>1$.
For the first case $0<s<1$,   the    discrete fractional Laplace operator $(-\Delta)^{s}$ acting on a function $u: V \rightarrow \mathbb{R}$ reads as
\begin{align*}
(-\Delta)^s u(x)=\frac{1}{\mu(x)}\sum_{y \in  V, \, y \neq  x}W_s(x, y) \left(u(x)-u(y)\right),
\end{align*}
where $W_s(x, y)$ is a symmetric positive function  for all  $x \neq y\in V$, see \eqref{Hs} below for details. Moreover, it is well known that the    discrete   Laplace $-\Delta$ acting on a function $u: V \rightarrow \mathbb{R}$  can be represented as
\begin{align*}
   - \Delta  u(x)=\frac{1}{\mu(x)}\sum_{y\sim x}w_{xy}(u(x)-u(y)).
   \end{align*}
Obviously, our definition  of $(-\Delta)^s$  is formally similar to $ - \Delta$. In addition,  we study some important properties.  For any $i=1,2,\cdots, n$ with $n=\sharp V$,  let $\lambda_i$  be the eigenvalue of   $-\Delta$, and $\phi_i$ is the corresponding orthonormal eigenfunction. We conclude that $\lambda_i$   to the power of $s$, written  $\lambda^s_i$, is the eigenvalue of  $(-\Delta)^{s}$ and $\phi_i$ remains the corresponding orthonormal eigenfunction,  namely
$$
( -\Delta)^s \phi_i =\lambda_i^s\phi_i , \ \forall i= 1, 2,\cdots,n.
$$
We then give an explicit expression of $(-\Delta)^s u $:
\begin{align*}
(-\Delta)^s u(x)= \sum_{i= {1}}^{ {n}}  \lambda_i^s\phi_i(x)  \langle u,\phi_i\rangle,
\end{align*}
  which is a property similar to $-\Delta$.
It should be remarked that in a general infinite graph, $(-\Delta)^s$ has no explicit expression.
Such an explicit expression is very important 
 for  both the theoretical analysis and computational procedures employed in problems involving  the fractional Laplacian.
 Furthermore, we study the asymptotic behavior as  ${s\rightarrow 1^-}$ and as ${s\rightarrow 0^+}$ for $(-\Delta)^{s}$: $  \lim_{s\rightarrow 1^-} (-\Delta)^{s} u=-\Delta u$  and
  $ \lim_{s\rightarrow 0^+} (-\Delta)^{s} u=u$
 if  $\int_V{u}d\mu=0$, which is similar to the property \eqref{a} in Euclidean  space.
We remark that the proof  in the graph setting is essentially different from that in the Euclidean setting.

As for the second case $s>1$,  we   rewrite $s =   \sigma+m$ with  $\sigma \in(0, 1)$ and $m \in \mathbb{N}_+$, and   define   $(-\Delta)^{s}$ in the distributional sense as
 \begin{align*}
\int_{V} \varphi (-\Delta)^s u  d\mu=  \int_{V}\nabla^s \varphi  \nabla^s u d\mu, \   \forall   \varphi \in C(V),
\end{align*}
where $\nabla^{s} u= \nabla^\sigma  \nabla^m u$.
We then give an explicit expression of  $(-\Delta)^{s}$:
 \begin{align}\label{s1}
(-\Delta)^{s}    u (x)=
\left\{\begin{aligned}
& -\Delta^{\frac{m-1}{2}}  \mathrm{div} (-\Delta)^\sigma  \nabla  \Delta^{\frac{m-1}{2}}   u(x) , & \text {if $m$ is odd,} \\
&  \Delta^{\frac{m}{2}}(-\Delta)^\sigma \Delta^{\frac{m}{2}}   u(x), & \text {if $m$ is eve.}\,
\end{aligned}\right.
 \end{align}
 Specially if $m=0$,   the above  expression   implies that  $(-\Delta)^{s}    u = (-\Delta)^\sigma u$ from the fact  $\Delta^0 u =u $.
To sum up, \eqref{s1} can be seen as a definition of $(-\Delta)^{s}$  for an arbitrary finite graph and any positive real number $s$.

 On the other hand,  we  pay attention to   the fractional Kazdan-Warner equation  on a connected finite graph $G = (V, E, \mu, w)$, namely
\begin{align}\label{KW}
(-\Delta)^{s} u=\kappa e^{u}-c \ \text { in } \ V,
\end{align}
where $c$ is a constant and $\kappa$ is some prescribed function on the graph.
 The solvability of  (\ref{KW}) will be discussed in two cases: $0<s<1$ and $s>1$.
For the first case $0<s<1$,  depending on the sign of $c$, we have the following
\begin{theorem}\label{T2}
 Let 
 $0<s<1$, $c$ be a  constant and $\kappa : V\rightarrow \mathbb{R}$ be
a function.  \\
(i) If $c>0$, then (\ref{KW}) is solvable if and only if $\kappa$ is positive somewhere.\\
(ii)  If $c=0$ and  $\kappa \not\equiv 0$, then  (\ref{KW}) is solvable if and only if both $\kappa$ changes sign and $\int_{V} \kappa d \mu<0$.\\
(iii)  If there exists  some $c_0<0$ such  that   (\ref{KW}) is solvable, then $\int_{V} \kappa d \mu<0$.
\end{theorem}
\begin{theorem}\label{T4}
 Let 
 $0<s<1$, $c<0$  and   $\kappa : V\rightarrow \mathbb{R}$ be a function satisfying    $\int_{V} \kappa d \mu<0$. \\
 (i) There exists a  constant $-\infty\leq c_{s,\kappa}<0$ depending only on $s$ and $\kappa$,  such that (\ref{KW}) is solvable for any $c_{s,\kappa}<c<0$, but   unsolvable for any $c<c_{s,\kappa}$.\\
(ii) $c_{s,\kappa}= -\infty$  if and only if  $\kappa(x) \leq 0$ for all $x \in V$.\\
(iii) If $c_{s,\kappa}\neq -\infty$, then  (\ref{KW}) is solvable when $ c =c_{s,\kappa}$.
 \end{theorem}
As for  $s>1$, we obtain the following theorem, which is weaker than the above theorems.
\begin{theorem}\label{T2-m}
 Let 
$s>1$, $c$ be a  constant and $\kappa : V\rightarrow \mathbb{R}$ be a function.  \\
(i) If $c>0$,  then  (\ref{KW}) has a solution if and only if $\kappa$ is positive somewhere.\\
(ii)  If $c=0$,  $\kappa$ changes sign and $\int_{V} \kappa d \mu<0$, then (\ref{KW})  is solvable.\\
(iii)  If $c<0$  and    $\kappa(x) <0$ for all $x \in V$, then (\ref{KW}) is solvable.
\end{theorem}
We will prove Theorems \ref{T2}--\ref{T2-m} using variational principles and the method of upper and lower solutions.  These methods were pioneered by Kazdan-Warner \cite{K-W-1, K-W-2} and later adapted to finite graphs by Grigor'yan-Lin-Yang \cite{A-Y-Y-2}.
Through the basic ideas of our proofs derived from them,   some technique difficulties caused by  $(-\Delta)^s$ must be overcome.   \\

The remaining part of this paper is organized as follows: In Section \ref{S2},  we make some preparations needed to state our main results.
In Section \ref{S3},  we define a discrete version of  the fractional Laplace  operator
and give some important properties.
In Section \ref{S6}, 
we study the solvability of the fractional   Kazdan-Warner equation.  For simplicity,  we do not distinguish between sequence and subsequence unless necessary. Moreover, we use the capital letter $C$ to denote some uniform constants that are independent of the special solutions
and not necessarily the same at each appearance.

\section{Preliminaries}\label{S2}

 We recall here some basic notions and important properties on graphs that will be used in the following discussions.
  Let   $V$ be a finite set of vertices,   $\sharp V$ be the number of all distinct vertices in $V$, and $E=\{[x,\,y]: x,y\in V  \text{ satisfy }  x\sim y\}$ be a finite set of edges, where $x\sim y$ means $x$ is connected to $y$ by an edge. For any vertex $x\in V$, we assign its measure $\mu:V\rightarrow \mathbb{R}^{+}$ with $x\mapsto \mu(x)$. For any edge $[x,\,y]\in E$, we denote its weight $w:E\rightarrow \mathbb{R}^{+}$ with  $[x,\,y]\mapsto  w_{xy}$, and
we always assume $w$ is symmetric, that is,  $w_{xy}=w_{yx}$ for any $x\sim y$.
 The quadruple $G = (V, E, \mu, w)$ refers to a finite graph.
Moreover, a graph is called connected if any two vertices can be connected via finite edges.
Throughout this paper, we always assume that $G = (V, E, \mu, {w})$ is a connected finite graph.

The  Laplace  operator acting on a function $u: V \rightarrow \mathbb{R}$ reads as
   \begin{align}\label{laplace}
   - \Delta  u(x)=\frac{1}{\mu(x)}\sum_{y\sim x}w_{xy}(u(x)-u(y)).
   \end{align}
In \cite{S-Y-Z-1},   the authors first formalized the gradient as a vector-valued function. Namely,
 the  gradient  form of  $u: V \rightarrow \mathbb{R}$ is written as
 \begin{align*}
 \nabla u(x)=\left(\sqrt{\frac{ w_{xy_1}}{2\mu(x)}}\left (u(x)-u(y_1)\right), \cdots, \sqrt{\frac{ w_{xy_{\ell_x}}}{2\mu(x)}}\left (u(x)- u(y_{\ell_x})\right)
\right),
\end{align*}
where  $y_i\sim x$ for all $i=1,2,\cdots,  \ell_x$ with  $\ell_x=\sharp \{y\in V:\,y\sim x\}$.
 The inner product of the gradient is assigned  as
 $$ \nabla u \nabla v   (x)
  =   \frac{1}{2{\mu(x)}}\sum_{y \sim x}   {w}_{xy} (u(x)-u(y))(v(x)-v(y))$$
and the length of $u$'s  gradient   is written  as $|\nabla u|  (x)=  \sqrt{\nabla  u  \nabla  u(x)}.  $
In addition, higher order derivatives can be found in   \cite{A-Y-Y-1, A-Y-Y-2, S-Y-Z-1}. Exactly,
 for any $m \in \mathbb{N}_+$,  the $m$-order  Laplacian operator  acting on $u: V \rightarrow \mathbb{R}$ is defined by
 $$\Delta^m u(x)=\Delta (\Delta^{m-1} u)(x)$$ and  the $m$-order gradient is
\begin{align}\label{Gradient-m}
 \nabla^m u=\left\{\begin{aligned}
& \nabla \Delta^{\frac{m-1}{2}} u , & \text {if $m$ is odd,}\ \,\\
& \Delta^{\frac{m}{2}} u , & \text {if $m$ is even.}
\end{aligned}\right.
\end{align}

We now introduce some concepts of integral and function spaces.  Let $C(V)$ be the space of all real-valued functions on the graph $G$.
	The integral of a function $u\in C(V)$ is defined by
	\begin{align*}
		\int_{V}ud\mu=\sum_{x\in V}u(x)\mu(x).
	\end{align*}
For any  $p\in [1,+\infty)$,    ${L^p(V)}$   denotes a linear space of all functions $u\in C(V)$ with finite norms
$$\|u\|_{L^p(V)}=\left(\sum_{x\in  V}|u(x)|^p\mu(x)\right)^{\frac{1}{p}},$$
and  ${L^\infty(V)}$ is a linear space of all bounded functions with the norm
  \begin{align*}
	 \|u\|_{L^\infty(V)}=\sup_{x\in  V}|u(x)|.
\end{align*}
 For brevity, we will use $\|\cdot\|_p$ and $\|\cdot\|_{\infty}$
 to denote $\|\cdot\|_{L^p(V)}$  and $\|\cdot\|_{L^\infty(V)}$ respectively.
Moreover,    the Sobolev space $W^{1, 2}(V)$ is
  \begin{align*}
W^{1, 2}(V)=\left\{u\in C(V): \int_V\left(|\nabla u|^2+u^2\right) d \mu<+\infty\right\},
\end{align*}
which is equipped with the norm
  \begin{align*}
\|u\|_{W^{1, 2}(V)}=\left(\|\nabla u\|_2^2+\|u\|_2^2\right)^{\frac{1}{2}}.
\end{align*}
In addition, for any $m \in \mathbb{N}_+$,   the set $$W^{m, 2}(V)=\left\{u\in C(V):  |\nabla^j u| \in L^2(V) \text{ for all } j=0,1, \cdots, m\right\}$$ is also a Sobolev space  with the norm
  \begin{align*}
\|u\|_{W^{m, 2}(V)}=\left(\sum_{j=0}^m \|\nabla^j u \|_2^2\right)^{\frac{1}{2}},
\end{align*}
where  $ \nabla^j u(x)$ is given as in \eqref{Gradient-m} for all $j=0,1, \cdots, m$.

 Moreover, on the   finite graph
 $G = (V, E, \mu, {w})$ without boundary constraint, the fundamental solution of the  discrete heat equation
  \begin{align*}\partial _t u(t,x)= \Delta u(t,x) \ \text{ on }\  [0,+\infty)\times V \end{align*}
is called  ``heat kernel" and is written as  $p(t, x, y) $.
 The heat kernel $p(t, x, y)$   can be represented as
 \begin{align}\label{heat-2}
 p(t,x,y)=\sum_{i=1}^{n} e^{-\lambda_it}\phi_i(x)\phi_i(y),
 \end{align}
 where  $n =\sharp V$,  $\lambda_i$ denotes the eigenvalue of   $-\Delta$, and
$\phi_i$   is the corresponding orthonormal eigenfunction.
Given a bounded initial condition $ u(0,x)=u_0(x)$  for any $x\in V$,  the above heat equation has a   uniquely bounded solution
\begin{align}\label{e2}
e^{t \Delta}u_0(x)
=\sum_{y\in V} p(t, x, y) u_0(y) \mu(y).\end{align}
It is not difficult to verify that the heat kernel $p(t, x, y)$ satisfies
 \begin{align}\label{heat-3}\underset{y\in V}\sum p(t,x,y)\mu(y)=1.\end{align}
The construction of heat kernel with Dirichlet or Neumann boundary conditions is carried out in  \cite{WA10, Wo08, Wo09} for infinite graphs or finite subgraphs,  with  $\mu\equiv1$ and   $w\equiv1$.
For the normal $\mu$ and $w$, we refer to \cite{Chung96, Keller-Lenz, HLLY}.
For our work,   we do not need additional boundary assumptions.
Although this is slightly different from the above situations, we can still find similar results from the book of Keller-Lenz-Wojciechowski  (\cite{Keller-book}, Chapter 0).

\section{Fractional Laplace operator}\label{S3}

This section is devoted to the definition of   $(-\Delta)^{s}$  for any real number $s>0$ on an arbitrary finite graph without any additional constraints.
Moreover, we give some important properties that will be used in the study of nonlinear problems involving the discrete operator $(-\Delta)^s$.
The discussion is presented in two cases: $0<s<1$ and $s>1$.

 \subsection{The case  $0<s<1$}

We will define the fractional Laplace operator $(-\Delta)^s $  on a connected finite graph $G = (V, E, \mu,  {w})$.  Specifically,   for any $0<s<1$,  it follows from \eqref{ee}, \eqref{e2} and \eqref{heat-3} that $(-\Delta)^{s}$ acting on a function $u\in C(V)$ can be characterized as
\begin{align*}
(-\Delta)^s u(x) \nonumber&=\frac{s}{\Gamma(1-s)} \int_0^{+\infty}\left(u(x)-e^{t \Delta}u(x)\right) t^{-1-s} d t\\
\nonumber&=\frac{s}{\Gamma(1-s)} \int_0^{+\infty}\left(1\cdot u(x)-\sum_{y \in  V}p(t, x, y) u(y)\mu(y)\right) t^{-1-s} d t\\
\nonumber&=\frac{s}{\Gamma(1-s)} \int_0^{+\infty}\sum_{y \in  V, \, y \neq  x}p(t, x, y) \mu(y)(u(x) - u(y) ) t^{-1-s} d t\\
 &=\frac{1}{\mu(x)}\sum_{y \in  V, \, y \neq  x}\left( \frac{s} {\Gamma(1-s)} \mu(x)\mu(y)\int_0^{+\infty}  p(t, x, y) t^{-1-s} d t\right) \left(u(x)-u(y)\right).
\end{align*}
Denote  a function
\begin{align}\label{Hs}
W_s(x, y)= \frac{s} {\Gamma(1-s)} \mu(x)\mu(y)\int_0^{+\infty}  p(t, x, y) t^{-1-s} d t ,\ \forall \ x \neq y \in  V.
\end{align}
Then for any $0<s<1$,   the fractional Laplace operator $(-\Delta)^s $ can be expressed on a finite graph as follows.
\begin{definition}\label{D1}
Let $G = (V, E, \mu,  {w})$ be a  connected finite graph and $0<s<1$.  The fractional Laplace  operator $(-\Delta)^{s}$ acting on a function $u\in C(V)$ reads as
\begin{align} \label{e-p}
(-\Delta)^s u(x)= \frac{1}{\mu(x)}\sum_{y \in  V, \, y \neq  x}W_s(x, y) \left(u(x)-u(y)\right),
\end{align}
where the   function $W_s(x, y)$, for any  $x \neq y\in V$, is defined as in \eqref{Hs}.
\end{definition}

 Let us now look at the properties of   $W_s(x, y)$, which is strictly related to $p(t, x, y)$.  From \eqref{heat-2}, the heat kernel $p(t, x, y)$ is represented  by   eigenfunctions and eigenfunctions of the    Laplace  operator $-\Delta$.
From  \cite{F-73, E-S-14},  the operator $-\Delta$ has an eigenvalue zero and the corresponding eigenfunction space has one dimension.      Without loss of generality, we require
  \begin{align}\label{lambda-2}
   0=\lambda_1<\lambda_2\leq \cdots\leq \lambda_{n}.
   \end{align}
     From this,
  we  calculate
  \begin{align}\label{phi1}
  \phi_1=\frac{1}{\sqrt{|V|}}\ \text{ and }\
 -\Delta  \phi_1(x)=0,
\end{align}
 where $|V|=\sum_{x\in V}\mu(x)$ denotes the volume of $V$,  and then
\begin{align}\label{lambda-0}
\sum_{x\in V}\phi_i(x)  \mu(x)
=0, \quad i= 2,\cdots,n.
    \end{align}
For any function $u\in C(V)$, there is
 \begin{align}\label{u}
    u=\sum_{i= {1}}^{ {n}}\langle u,\phi_i\rangle \phi_i,
    \end{align}
and thus    the characteristic function
 \begin{align}\label{delta}
 {1_x}(y)=\sum_{i= {1}}^{ {n}}\left\langle {1_x} ,\phi_i\right\rangle\phi_i(y)=
    \sum_{i= {1}}^{ {n}}\phi_i(x)\phi_i(y){\mu(x)}.
    \end{align}
Reviewing    \eqref{heat-2} and \eqref{lambda-2}, we derive that
 \begin{align}\label{h}
\nonumber p(t,x,y)  =&  \sum_{i=1}^{n} ( e^{-\lambda_it}-1)\phi_i(x)\phi_i(y)+ \sum_{i=1}^{n}\phi_i(x)\phi_i(y)\\
 =&\sum_{i=2}^{n}( e^{-\lambda_it}-1)\phi_i(x)\phi_i(y)+  \frac{1_x}{\mu(x)}(y).
 \end{align}
Therefore,  $W_s(x, y)$  can be explicitly expressed in terms of eigenfunctions and eigenfunctions of the discrete  Laplace operator $-\Delta$ as follows.
 \begin{proposition}\label{Pw}
The function $W_s(x, y)$ in \eqref{Hs} can be expressed as
\begin{align}\label{H}  W_s(x, y) =  -\mu(x)\mu(y)\sum_{i=1}^{n}\lambda_i^s\phi_i(x)\phi_i(y),\ \forall\, x \neq y\in V, \end{align}
where $n=\sharp V$,   $\lambda_i$ is the eigenvalue of the Laplace  operator $-\Delta$, and
$\phi_i$ is the corresponding orthonormal eigenfunction. Moreover,  the function $ W_s(x, y)$ satisfies
\begin{align}\label{H2}
0<W_s(x, y)=W_s(y,x)<+\infty,\ \forall\, x \neq y\in V. \end{align}
\end{proposition}

\begin{proof}
Inserting \eqref{h} into \eqref{Hs},  we conclude
 for any $ x \neq y\in V $,
 \begin{align}\label{W}
   W_s(x, y) 
        =&   \frac{s} {\Gamma(1-s)} \mu(x)\mu(y) \int_0^{+\infty}  \left( \sum_{i=2}^{n}( e^{-\lambda_it}-1)\phi_i(x)\phi_i(y) \right) t^{-1-s} d t    \nonumber\\
        =&   \frac{s} {\Gamma(1-s)} \mu(x)\mu(y)  \sum_{i=2}^{n}\phi_i(x)\phi_i(y)   \int_0^{+\infty}(e^{-\lambda_it}- 1) t^{-1-s} d t  \nonumber  \\
        \nonumber =&   \frac{s} {\Gamma(1-s)}   \mu(x)\mu(y)\sum_{i=2}^{n}\lambda_i^s\phi_i(x)\phi_i(y)   \int_0^{+\infty}(e^{-t}- 1) t^{-1-s} d t\\
         =&-\mu(x)\mu(y)\sum_{i=2}^{n}\lambda_i^s\phi_i(x)\phi_i(y) .
\end{align}
The last equality is due to
      \begin{align*}
     \int_0^{+\infty}  (e^{-t}- 1)  t^{-1-s} d t =  -\frac{1}{s}\int_0^{+\infty}  (e^{-t}- 1)    d t^{-s}
 =  -\frac{\Gamma(1-s)}{s} \end{align*}
from integration by parts.  Then \eqref{H} follows from \eqref{W}.
Moreover,
 \eqref{H2} is obvious from  \eqref{Hs} and  \eqref{H}.
\end{proof}

Moreover, we claim that $\phi_i$ remains an orthonormal eigenfunction corresponding to the eigenvalue  $\lambda^s_i$  of the fractional Laplace operator $(-\Delta)^{s}$,  namely the following
\begin{proposition}\label{T}
Let $0<s<1$.
For any $x\in V$, there holds
 \begin{align*}
( -\Delta)^s \phi_i =\lambda_i^s\phi_i , \ \forall i= 1, 2,\cdots,n.
    \end{align*}
\end{proposition}
\begin{proof}  For any $x\in V$ and $ i =2$, $\cdots$, $n$,   it follows from \eqref{e-p}--\eqref{phi1} and \eqref{H}   that
    \begin{align}\label{8}
(-\Delta)^s \phi_i(x)
= &-\sum_{y \in  V, \, y \neq  x}  \sum_{j=2}^{n}\lambda_j^s\phi_j(x)\phi_j(y) \mu(y)\left(\phi_i(x)-\phi_i(y)\right)\nonumber\\
= & \sum_{j=2}^{n}  \lambda_j^s\phi_j(x)\left(-\phi_i(x)\sum_{y \in  V, \, y \neq  x}\phi_j(y) \mu(y)+\sum_{y \in  V, \, y \neq  x}\phi_i(y)\phi_j(y) \mu(y)\right).
\end{align}
Since \eqref{lambda-0},  for any $ j =2$, $\cdots$, $n$, we obtain
    \begin{align*}
\sum_{y \in  V, \, y \neq  x}\phi_j(y) \mu(y)=   \sum_{y \in  V}\phi_j(y) \mu(y)-\phi_j(x) \mu(x) =-\phi_j(x) \mu(x).
\end{align*}
 Inserting it into   \eqref{8}, we get
     \begin{align}\label{2}
\nonumber(-\Delta)^s \phi_i(x)
= &  \sum_{j=2}^{n}  \lambda_j^s\phi_j(x)\left( \phi_i(x) \phi_j(x) \mu(x)+\sum_{y \in  V, \, y \neq  x}\phi_i(y)\phi_j(y) \mu(y)\right)\\
= &  \sum_{j=2}^{n}  \lambda_j^s\phi_j(x) \sum_{y \in  V}\phi_i(y)\phi_j(y) \mu(y).
\end{align}
Since the  eigenfunctions $\phi_i$ ($1\leq i\leq n$) are  orthonormal, there holds for any $ i,j= 1,2,\cdots,n,$
    \begin{align*}
    \langle {\phi}_i, {\phi}_j\rangle=\sum_{y\in V}\phi_i(y)\phi_j(y) \mu(y)
    =\left\{\begin{aligned}
    &1,&{\rm if }\ \ i=j,\\
    &0,&{\rm if }\ \ i\not=j,
    \end{aligned}\right.
    \end{align*}
which together with \eqref{2} leads to this proposition.
\end{proof}

Therefore,  $(-\Delta)^s u $ can also be explicitly expressed in terms of eigenfunctions and eigenfunctions of  $-\Delta$ as follows.
\begin{proposition} \label{P1}
For any   $u\in C(V)$,    $(-\Delta)^{s}u$ in \eqref{e-p}   can be expressed as
\begin{align}\label{e-p-2}
(-\Delta)^s u(x)= \sum_{i= {1}}^{ {n}}  \lambda_i^s\phi_i(x)  \langle u,\phi_i\rangle,
\end{align}
where $n=\sharp V$,   $\lambda_i$ is the eigenvalue of the Laplace  operator $-\Delta$, and $\phi_i$ is the corresponding orthonormal eigenfunction.
\end{proposition}
\begin{proof}
It follows from   \eqref{e-p}   and \eqref{H}  that
\begin{align}\label{1-1}
(-\Delta)^s u(x)
= &\frac{1}{\mu(x)}\sum_{y \in  V, \, y \neq  x}W_s(x, y) \left(u(x)-u(y)\right)\nonumber\\
= & \sum_{y \in  V, \, y \neq  x}- \mu(y)\sum_{i=1}^{n}\lambda_i^s\phi_i(x)\phi_i(y) \left(u(x)-u(y)\right)\nonumber\\
= & \sum_{i=1}^{n} \lambda_i^s\phi_i(x)  \sum_{y \in  V } \phi_i(y)  u(y) \mu(y)-u(x) \sum_{y \in  V }\sum_{i=1}^{n}  \lambda_i^s\phi_i(x) \phi_i(y) \mu(y).
\end{align}
In view of  \eqref{delta}, we have $  {1_y}(x) =    \sum_{i= {1}}^{ {n}}\phi_i(x)\phi_i(y) {\mu(y)}.$ Then  from  Proposition \ref{T} and the  linearity of $(-\Delta)^s$, there holds
\begin{align*}
\sum_{y \in  V }\sum_{i=1}^{n}  \lambda_i^s\phi_i(x) \phi_i(y) \mu(y)=&\sum_{y \in  V }\sum_{i=1}^{n}  ( -\Delta)_x^s \phi_i(x) \phi_i(y) \mu(y)\\
= & ( -\Delta)_x^s \sum_{y \in  V }\sum_{i=1}^{n}  \phi_i(x) \phi_i(y) \mu(y)\\
= & ( -\Delta)_x^s \sum_{y \in  V }{1_y}(x) \\
= & 0.
\end{align*}
Inserting it into \eqref{1-1},  we conclude \eqref{e-p-2}.
\end{proof}

Moreover, an important property of $(-\Delta)^{s}$ is its asymptotic behavior as ${s\rightarrow 1^-}$ and ${s\rightarrow 0^+}$.
Similar to  the Euclidean fractional Laplace operator, we have the following
 \begin{proposition}\label{sT1}
Let $0<s<1$. For any $u\in C(V)$, there hold:\\
    (i) $ \underset{s\rightarrow 1^-}\lim  (-\Delta)^{s} u=-\Delta u$;\\
(ii) $\underset{s\rightarrow 0^+} \lim (-\Delta)^{s} u=u$
 if $u$ satisfying $\int_V{u}d\mu=0$.
\end{proposition}
\begin{proof}
 On the other hand,  it follows from  Proposition  \ref{P1} and the  linearity of $(-\Delta)^s$ that
\begin{align*}
\lim_{s\rightarrow 1^-}(-\Delta)^s u(x)  = &\lim_{s\rightarrow 1^-}  \sum_{i= {1}}^{ {n}}  \lambda_i^s\phi_i(x)  \langle u,\phi_i\rangle\\
= &  \sum_{i= {1}}^{ {n}}  \lambda_i \phi_i(x)  \langle u,\phi_i\rangle\\
= &    -\Delta_x\sum_{i= {1}}^{ {n}} \phi_i(x)  \langle u,\phi_i\rangle,
\end{align*}
which together with \eqref{u} leads to ($i$).
On the other hand,   from   \eqref{lambda-2} and \eqref{phi1}, we get
\begin{align*}
\lim_{s\rightarrow 0^+}(-\Delta)^s u(x)  
= &\lim_{s\rightarrow 0^+} \sum_{i= {2}}^{ {n}}  \lambda_i^s\phi_i(x)  \langle u,\phi_i\rangle+  0\\
= & \sum_{i= {2}}^{ {n}}  \phi_i(x)  \langle u,\phi_i\rangle+ \phi_1(x)  \langle u,\phi_1\rangle- \phi_1(x)  \langle u,\phi_1\rangle\\
= & \sum_{i= 1}^{ {n}}  \phi_i(x)  \langle u,\phi_i\rangle - \phi_1(x)  \langle u,\phi_1\rangle\\
= & u(x) - \bar{u}.
\end{align*}
If $u$ satisfying $\int_V{u}d\mu=0$,   ($ii$) follows.
\end{proof}

\begin{remark}
 Building upon the previous analysis, we obtain that the fractional Laplace operator admits a matrix representation.
Proposition \ref{sT1}  is then concluded  through standard techniques of matrix theory,  and a complete verification is provided  in the Appendix.
\end{remark}

In addition,  noting that $V$ is finite, for any fixed $x\in V$, we  write $V=\{x, y_1, y_2,\cdots,y_{n-1}\}$ with $n=\sharp V$.
Then  the  fractional gradient   form of  $u\in C(V)$ reads as
\begin{align}\label{Gradient-s}
 \nabla^s u(x)=\left(\sqrt{\frac{ {W_s(x, y_1)}}{2\mu(x)}}\left (u(x)-u(y_1)\right), \cdots, \sqrt{\frac{ {W_s(x, y_{n-1})}}{2 \mu(x)}} \left (u(x)-u(y_{n-1})\right)
\right).
\end{align}
The inner product of the  fractional gradient is assigned  as
 $$ \nabla^su \nabla^sv   (x)
  =   \frac{1}{2\mu(x)}\sum_{y \in  V, \, y \neq x} W_s(x, y) (u(x)-u(y))(v(x)-v(y))$$
and the length of $u$'s fractional gradient   is written  as
\begin{align}\label{N}
|\nabla^su|  (x)=  \sqrt{\nabla^s u  \nabla^s u(x)}
  =\left( \frac{1}{2\mu(x)}\sum_{y \in  V, \, y \neq x} W_s(x, y) (u(x)-u(y))^2 \right)^{\frac{1}{2}}.
  \end{align}

  We next focus on the formula of integration by parts, which is the key to using the calculus of variations.
 \begin{proposition}\label{Lpart}
Let  $0<s<1$.
For any $u$, $ v\in C(V)$, there holds
\begin{align}\label{part}
\int_{ V} v(-\Delta)^s u  d\mu=  \int_{ V}  \nabla^s u \nabla^s v d\mu =\int_{ V}u(-\Delta)^{s} v  d\mu,
\end{align}
where  $(-\Delta)^s$ and   $\nabla^s $ are defined by
\eqref{e-p} and \eqref{Gradient-s} respectively.
\end{proposition}
\begin{proof}
By a direct computation, we derive that for any  $u, v\in C(V)$,
\begin{align}\label{1}
\nonumber\int_{ V} v(-\Delta)^s u  d\mu =& \sum_{x\in V}v(x) \sum_{{y \in  V, \,y \neq x}}W_s(x, y) (u(x)-u(y))\\
\nonumber=& \sum_{x\in V}\sum_{{y \in  V, \,y \neq x}}(v(x)-v(y)) (u(x)-u(y))W_s(x, y) \\
\nonumber&+\sum_{x\in V}\sum_{{y \in  V, \,y \neq x}}  v(y)  (u(x)-u(y))W_s(x, y) \\
=&2\int_{ V} \nabla^s u \nabla^s vd\mu+\sum_{x\in V}\sum_{{y \in  V, \,y \neq x}}v(y)  (u(x)-u(y))W_s(x, y) .
\end{align}
It follows from   $G$ is a  connected finite graph and \eqref{H2} that
\begin{align*}
\sum_{x\in V}\sum_{{y \in  V, \,y \neq x}}v(y)  (u(x)-u(y))W_s(x, y)
&=\sum_{y\in V}\sum_{x \in V,\,x \neq y}v(y)  (u(x)-u(y))W_s(x, y) \\
&=\sum_{x\in V}\sum_{{y \in  V, \,y \neq x}}v(x)  (u(y)-u(x))W_s(x, y) \\
&=-\int_{ V} v(-\Delta)^s u  d\mu.\end{align*}
Inserting it into   \eqref{1}, we conclude this result.
\end{proof}

Finally in this subsection, we give some important properties below.

\begin{proposition}\label{P2}
Let $0<s<1$. For any $u \in C(V)$, the following statements hold:\\
(i)   if  $ \nabla^{s} u(x_0)=\mathbf{0}$ for some vertex $x_0\in V$, where $\mathbf{0}$ denotes the zero vector,  then $u$ is a constant value function;\\
(ii)  if $ (-\Delta)^{s} u(x)\leq0$ for all $x\in V$, and $u(x_0)=\max_Vu$ for some  vertex $x_0\in V$, then $u$ is a constant value function;\\
(iii)   if $ (-\Delta)^{s} u\equiv0$, then $u$ is a constant value function;\\
 (iv)  $(-\Delta)^s (ku +\ell v ) =k (-\Delta)^s  u +\ell (-\Delta)^s v  $ for any constants $k$ and $\ell$;\\
(v) $(-\Delta)^s (u v ) =u (-\Delta)^s   v +v (-\Delta)^s  u +2 \nabla^su \nabla^sv    $ for any $u,\, v \in C(V)$.
\end{proposition}

\begin{proof}
If  $ \nabla^{s} u(x_0)=\mathbf{0}$ for some vertex $x_0\in V$,  then   from \eqref{H2} and  \eqref{Gradient-s}, we obtain  $u(y)=u(x_0)$ for any $     y \neq  x _0$, i.e.,  $u$ is a constant value function, and then ($i$) follows.

Furthermore,
it follows from  \eqref{e-p}  that
\begin{align*}
(-\Delta)^s u(x_0)
= \frac{1}{\mu(x_0)}\sum_{y \in  V, \, y \neq  x_0} {W_s(x_0, y) } \left(u(x_0)-u(y)\right)\geq0.
\end{align*}
Since  $ (-\Delta)^{s} u(x)\leq0$ for all $x\in V$, then $(-\Delta)^s u(x_0)=0$. Therefore,  $u$ is a constant value function, and then there is ($ii$).

 Moreover,  if $ (-\Delta)^{s} u(x)=0$ for all $x\in V$, then  it follows  from Proposition \ref{Lpart} that
\begin{align*}
0=\int_{ V} u(-\Delta)^s u  d\mu=  \int_{ V}  |\nabla^s u|^2d\mu, 
\end{align*}
and then $  |\nabla^s u|(x)=0$ for all $x\in V$. Because of \eqref{N}, we get  $u$ is a constant value function, and thus ($iii$) holds.

 In addition, ($iv$) follows from the definition of $(-\Delta)^s$ in \eqref{e-p}. And for any $x\in V$, a straightforward calculation leads to that
 \begin{align*}
 (-\Delta)^s( u v)(x)&=\frac{1}{\mu(x)}\sum_{y \in  V, \, y \neq  x}W_s(x, y) \Big(u(x)v(x)-u(x)v(y)+u(x)v(y)-u(y)v(y)\Big) \\
 &=u(x) (-\Delta)^s   v(x)+\frac{1}{\mu(x)}\sum_{y \in  V, \, y \neq  x}W_s(x, y) \left(u(x) -u(y)\right)v(y)  \\
  &=u(x) (-\Delta)^s   v(x)+v(x) (-\Delta)^s   u(x)+ 2 \nabla^su \nabla^sv(x).
\end{align*}
Then ($v$) holds, and ends the proof.
\end{proof}

\subsection{The case  $s > 1$}
If the fractional exponent $s > 1$, we can rewrite $s =   \sigma+m$ with  $\sigma \in(0, 1)$ and $m \in \mathbb{N}_+$. We define the fractional gradient form of  $u\in C(V)$  by
 \begin{align}\label{Gradient-m-s}
 \nabla^{s} u= \nabla^\sigma  \nabla^m u=\left\{\begin{aligned}& \nabla^\sigma  \nabla \Delta^{\frac{m-1}{2}} u , & \text {if $m$ is odd,}\ \,\\& \nabla^\sigma \Delta^{\frac{m}{2}} u, & \text {if $m$ is even.}\end{aligned}\right.
 \end{align}
Here $  \nabla^\sigma$   acting on a function is given by \eqref{Gradient-s}  and   acting on a $k$-dimensional vector-valued function $\mathbf{f}  =\left(f_1, f_2,\cdots, f_k\right): V \rightarrow \mathbb{R}^k$  is also a vector-valued function
  \begin{align*}
\nabla^{\sigma} \mathbf{f}(x)=\Big(\nabla^{\sigma}f_1(x), \nabla^{\sigma}f_2(x), \cdots,\nabla^{\sigma} f_k(x)\Big).
\end{align*}
Similarly,   $(-\Delta)^\sigma$   acting on a function is shown as in \eqref{e-p}   and   acting on a $k$-dimensional vector-valued functionn $\mathbf{f}  =\left(f_1, f_2,\cdots, f_k\right): V \rightarrow \mathbb{R}^k$  is also a vector-valued function
\begin{align*}
(-\Delta)^\sigma \mathbf{f} =\Big((-\Delta)^\sigma f_1 ,  (-\Delta)^\sigma f_2 , \cdots,(-\Delta)^\sigma f_k  \Big).
\end{align*}
Therefore, it follows from \eqref{part} that
\begin{align}\label{part-vector}
\int_{ V} \mathbf{f}\cdot (-\Delta)^\sigma  \mathbf{g}  d\mu=  \int_{ V}  \nabla^\sigma  \mathbf{f} \cdot \nabla^\sigma  \mathbf{g} d\mu =\int_{V}  \mathbf{g}\cdot (-\Delta)^\sigma \mathbf{f}  d\mu,
\end{align}
for any two vector-valued functions $\mathbf{f} $ and $\mathbf{g} $ of the same dimension.
Moreover, we define the divergence ``$\mathrm{div}$"  acting on      $\mathbf{f}$ by
\begin{align}\label{div}
\int_V (\mathrm{div} \mathbf{f})\, \varphi \,d \mu=-\int_V \mathbf{f}\cdot  \nabla \varphi  \,d \mu, \quad \forall \varphi \in C(V).
\end{align}
For any fixed $x \in V$, by taking $ \varphi$ as the Dirac function $1_x/{\mu(x)}$, we calculate
\begin{align}\label{div-B}
 \mathrm{div}   \mathbf{f} (x)=-\frac{1}{\mu(x)}\sum_{y\in V} \mu(y) \nabla_y 1_x  (y)   \cdot \mathbf{f}(y) ,
\end{align}
and then it is not difficult to check
\begin{align*}
\Delta u(x)=\mathrm{div}  \nabla u(x), \ \forall u\in C(V),
\end{align*}
 where $ \Delta  $ is defined as in \eqref{laplace}.

 We now  define the fractional Laplace  operator $(-\Delta)^{s}$ in the distributional sense  by
  \begin{align}\label{part-B}
\int_{V} \varphi (-\Delta)^s u  d\mu=  \int_{V}\nabla^s \varphi  \nabla^s u d\mu, \   \forall   \varphi \in C(V).
\end{align}
Next, we give an explicit expression for    $(-\Delta)^{s}$.
\begin{proposition}\label{D2}
 The fractional Laplace    operator $(-\Delta)^{s}$ acting on a function $u\in C(V)$ reads as
 \begin{align}\label{Delta-m-s}
(-\Delta)^{s}    u =
\left\{\begin{aligned}
& -\Delta^{\frac{m-1}{2}}  \mathrm{div} (-\Delta)^\sigma  \nabla  \Delta^{\frac{m-1}{2}}   u , & \text {if $m$ is odd,}\ \ \\
&  \Delta^{\frac{m}{2}}(-\Delta)^\sigma \Delta^{\frac{m}{2}}   u, & \text {if $m$ is even,}\,
\end{aligned}\right.
 \end{align}
where  $(-\Delta)^\sigma$ and ``$\mathrm{div}$"  are defined as in \eqref{e-p} and \eqref{div-B} respectively.
\end{proposition}
\begin{proof}
Let $u$ and $ \varphi \in C(V)$.
If $m$ is  odd,  it follows from  \eqref{Gradient-m-s}, \eqref{part-vector} and \eqref{div}  that
\begin{align*}
 \int_{V}\nabla^s \varphi  \nabla^s u d\mu&=
 \int_{ V}\nabla^{\sigma}  \nabla\Delta^{\frac{m-1}{2}}   \varphi \cdot \nabla^{\sigma} \nabla \Delta^{\frac{m-1}{2}}  u   d\mu\\
&=\int_{ V}   \nabla \Delta^{\frac{m-1}{2}}  \varphi  \cdot   (-\Delta)^{\sigma} \nabla\Delta^{\frac{m-1}{2}}   u  \, d\mu\\
&=-\int_{ V}  \Delta^{\frac{m-1}{2}}  \varphi \left(
 \mathrm{div} (-\Delta)^{\sigma}  \nabla  \Delta^{\frac{m-1}{2}}   u\right)    d\mu\\
  &=-\int_{ V}    \varphi \left(\Delta^{\frac{m-1}{2}}
 \mathrm{div} (-\Delta)^{\sigma}  \nabla  \Delta^{\frac{m-1}{2}}   u\right)    d\mu.
\end{align*}
Similarly, if $m$ is even, there holds
\begin{align*}
 \int_{V}\nabla^s \varphi  \nabla^s u d\mu&=
 \int_{ V}\nabla^{\sigma}      \Delta^{\frac{m}{2}}\varphi  \cdot\nabla^{\sigma}    \Delta^{\frac{m}{2}}   u   d\mu\\
&=\int_{ V} \Delta^{\frac{m}{2}}\varphi
\left( (-\Delta)^{\sigma}   \Delta^{\frac{m}{2}}   u \right)  d\mu\\
  &= \int_{ V} \varphi\left(\Delta^{\frac{m}{2}}(-\Delta)^{\sigma} \Delta^{\frac{m}{2}}   u \right)        d\mu.
\end{align*}
Then we conclude \eqref{Delta-m-s}  from \eqref{part-B}.
 \end{proof}
Specially if $m=0$,   the above definition  implies $(-\Delta)^{s}    u = (-\Delta)^\sigma u$ from the fact  $\Delta^0 u =u $. Then for  any fraction $ s>0$,
   $(-\Delta)^{s}$ acting on $u\in C(V)$
 can be expressed as in \eqref{Delta-m-s}.

\section{Fractional Kazdan-Warner equation}\label{S6}

In this section, we will prove Theorems \ref{T2}--\ref{T2-m} respectively. Before proving  these theorems
 formally, we present three embeddings in the next subsection, which help us to obtain satisfying characterizations of the solvability of the fractional  Kazdan-Warner equation on finite graphs.

\subsection{Embeddings} \label{S5}

 For any  $0<s<1$, the fractional Sobolev space $W^{s,2}(V)$ is defined by
 \begin{align*}
W^{s,2}(V) =\left\{u \in C(V) :\int_{V}\left(|\nabla^su|^2 +u^2\right)  d\mu<+\infty\right\},
\end{align*}
 equipped with the associated norm $\|u\|_{W^{s,2}(V)} = (\int_{V} (|\nabla^su|^2 +u^2 )  d\mu )^{1/{2}}$. 
For any $s >1$,   we   rewrite $s =   \sigma+m$ with  $\sigma \in(0, 1)$ and $m \in \mathbb{N}_+$, and define the  fractional Sobolev space
 \begin{align*}W^{s,2}(V)=\left\{u\in C(V):  |\nabla^{\sigma+j} u| \in L^2(V)\, \text{ and }\, |\nabla^j u| \in  L^2(V)  \text{ for all } j=0,1,  \cdots, m\right\}\end{align*} equipped with  the  norm
 \begin{align*}
 \|u\|_{W^{s,2}(V)}
 =\left(\sum_{j=0}^m  \|\nabla^{\sigma+j} u  \|_2^2+\sum_{j=0}^m \|\nabla^j u \|_2^2 \right)^{\frac{1}{2}},
 \end{align*}
 where   $ \nabla^j u(x)$ and $ \nabla^{\sigma+j} u(x)$  are given as in \eqref{Gradient-m} and \eqref{Gradient-m-s}  for all $j=0,1, \cdots, m$.

The Sobolev embeddings of fractional Sobolev space in Euclidean space are elaborated and proved in \cite{A}, and we highly recommend it to interested readers.  However,
these become easier on finite graphs.
Since $G $ is a connected finite graph, it is obvious that $W^{s,2}(V)$ is exactly the set of all functions on $G$
for any $s>0$, then  $W^{s,2}(V)$  is a finite-dimensional linear normed space, and is also a Banach space.  This implies the following Sobolev embedding.

 \begin{lemma}\label{Lse}
For any $s>0$, the  space $W^{s,2}(V)$ is pre-compact.
Namely, if $\{u_n\}$ is bounded in $W^{s,2}(V)$, then there exists some $u \in W^{s,2}(V)$ such that  $u_n \rightarrow u$ in $W^{s,2}(V)$.
\end{lemma}

It is well known that any two norms on a finite-dimensional space must be equivalent.  Since $G$ is a finite graph,
 we obtain the following Poincar\'{e} inequality.
 \begin{lemma}\label{Lpoincare}
For any  $s>0$ and any function  $u\in C(V)$ with $\int_{V} u d \mu=0$,
there exists some constant $C>0$ depending only on $G$ such that
  \begin{align*}
\int_{V}    u ^2 d \mu \leq C \int_{V} |\nabla^s u|^2  d \mu.
\end{align*}
\end{lemma}

In addition, we introduce the following Trudinger-Moser embedding.

\begin{lemma}\label{LTM}
For any  $s>0$ and any $\alpha \in \mathbb{R}$, there exists a positive constant $C$
depending only on $\alpha$ and $G$ such that for all functions $u\in C(V)$ with $\int_{V}|\nabla^s u|^{2} d \mu \leq 1$ and $\int_{V} u d \mu=$
0, there holds
  \begin{align*}
\int_{V} e^{\alpha u^{2}} d \mu \leq C.
\end{align*}
\end{lemma}

\begin{proof}
Since the case $\alpha \leq 0$ is trivial,  we will focus on the case of $\alpha>0$. For any function $u$ satisfying
$\int_{V}|\nabla^s u|^{2} d \mu \leq 1$ and $\int_{V} u d \mu=0$, it follows from Lemma \ref{Lpoincare} that
 $
\int_{V}  u ^{2} d \mu
\leq C
$
for some constant $C$ depending only on  $G$.
Denote $\mu_{\min }=\min _{x \in V} \mu(x)$.
Then the above inequality leads to
$ u^{2}(x)
\leq  {C}/{\mu_{\min }}
$
for any $x\in V$.
Hence
$
\int_{V} e^{\alpha u^{2}} d \mu \leq e^{{ \alpha C} /{ \mu_{\min }}} |V|,
$
which gives the desired result.
\end{proof}

\subsection{Proof of Theorem \ref{T2}}\label{sub4.1}
In this subsection,  we prove   Theorem \ref{T2} ($i$)-($iii$) respectively. At first, we consider the case of $c>0$.

\textit{\textbf{Proof of Theorem \ref{T2} ($i$).}}
{ \textit{Sufficiency.}}
Suppose $c>0$ and $u$ is a solution of (\ref{KW}). It follows from $\int_V (-\Delta)^{s} u d \mu=0$ that
$
\int_V  \kappa e^u d \mu=c |V|>0 .
$
Then $\kappa$ must be positive somewhere on $V$.

{ \textit{Necessity.}}
  Define a function space
\begin{align}\label{lambda-1}
\Lambda_1=\left\{ {u \in}  {W^{s,2}(V)}: \int_V  \kappa e^u d \mu=c|V|\right\} .
\end{align}
Suppose $\kappa \left(x_0\right)>0$ for some $x_0 \in V$. Then
it follows from (\cite{A-Y-Y-2}, Section 5) that  $\Lambda_1 $ is not empty.
Define
\begin{align}\label{J1}
 J_1(u)=\frac{1}{2} \int_V|\nabla^s u|^{2} d \mu+c \int_V u d \mu,\ \forall\ u\in \Lambda_1,
\end{align}
and $\theta_1=\inf _{ {u \in} \Lambda_1}  J_1(u)$.
To ensure the reasonableness of  $\theta_1$,
we next  prove that $J_1$ has a lower bound on $\Lambda_1$.
Write $v=u-\bar{u}$, then $\bar{v}=0$ and $\|\nabla^s v\|_2=\|\nabla^s u\|_{2}$. Hence for any $u\in \Lambda_1$,  we have
$\int_V  \kappa e^v d \mu= c |V| e^{-\bar{u}}>0 $,
and thus
\begin{align}\label{e1}
\int_V u d \mu=  |V| \left(\ln (c |V|) -  \ln \int_V  \kappa e^v d \mu\right).
\end{align}
Let $\tilde{v}=v /\| {\nabla^s v}\|_{2}$, then $\int_V \tilde{v} d \mu=0$ and $\|\nabla^s \tilde{v}\|_{2}=1$.
From Lemma \ref{Lpoincare}, there holds $\|\tilde{v}\|_{2} \leq C$.
For any $\alpha \in \mathbb{R}$,
it follows from Lemma \ref{LTM} that $\int_V e^{\alpha  \tilde{v} ^{2}} d \mu \leq C$.  For any $\epsilon>0$,
this together with   Young's  inequality $a b \leq  \epsilon a^{2} +  b^2/(4 \epsilon)$ implies that
\begin{align*}
\int_V e^v d \mu
 \leq \int_V e^{ \epsilon  \| {\nabla^s v}\|_{2}^{2}+\frac{1}{4 \epsilon} \tilde{v} ^2} d \mu
 \leq C  e^{\epsilon  \| {\nabla^s u}\|_{2}^{2}}.
\end{align*}
In view of (\ref{e1}),  the above inequality leads to
$
\int_V u d \mu 
  \geq -|V|  (    {\epsilon  \| {\nabla^s u}\|_{2}^{2}} +   C ).
$
Choosing $\epsilon=({2 c |V|})^{-1}$,  we obtain that
$J_1$ has a lower bound on the set $\Lambda_1$.

Take a sequence of functions $\left\{u_i\right\} \subset \Lambda_1$ such that $J_1(u_i) \rightarrow \theta_1$ as ${ i \rightarrow +\infty}$.  Moreover, the equality
  \begin{align*}
\int_V u_i d \mu=\frac{1}{c} J_1(u_i)-\frac{1}{2c} \int_V\left| {\nabla^s u}_i\right|^{2} d \mu
\end{align*}
implies that $\left\{\bar{u}_i\right\}$ is a bounded sequence.
Write $v_i=u_i-\bar{u}_i$, then $\bar{v}_i=0$ and $\|\nabla^s v_i\|_{2}=\|\nabla^s u_i\|_{2}$.
Hence   it follows from Lemma \ref{Lpoincare} that $v_i$ is bounded in $ {W^{s,2}(V)}$. Then $ u_i=v_i+\bar{u}_i$ is also bounded in $ {W^{s,2}(V)}$.
By the Sobolev embedding  in Lemma \ref{Lse},   there exists some $u_0 \in W^{s,2}(V)$ such that  $ {u_i} \rightarrow  u_0$ in $ {W^{s,2}(V)}$ as ${i \rightarrow +\infty}$.
It is easy to see that   $u_0$ is a minimize of $ J_1(u)$ on $\Lambda_1$, namely $u_0 \in \Lambda_1$ and $ J_1(u_0)=\theta_1$.  Moreover,   we derive the Euler-Lagrange equation of the minimizer $u_0$ as
 $(-\Delta)^{s} u_0=  \kappa e^{u_0}-c,$
which implies that $u_0$ is a solution of  (\ref{KW}). $\hfill\Box$\\

If $c = 0$, then (\ref{KW}) is reduced to
\begin{equation}\label{KW-1}(-\Delta)^{s} u=\kappa e^{u} \ \text { in }\ V.\end{equation}

\textit{\textbf{Proof of Theorem \ref{T2} ($ii$).}}
{\textit{Sufficiency.}}  If there is a solution $u$ to (\ref{KW-1}), then $u$ cannot be a constant function.   Otherwise,   it will conflict with the condition of  $\kappa \not\equiv 0$.  Integration of  (\ref{KW-1}) yields
  \begin{align*}
\int_V \kappa e^{u} d \mu=\int_V (-\Delta)^{s} u d \mu=0 .
\end{align*}
This together with $\kappa \not \equiv 0$ implies that $\kappa$ must change sign. On the other hand,  it is easy to get  $\kappa=e^{-u} (-\Delta)^{s} u$ from (\ref{KW-1}).  Integration by parts in  \eqref{part} gives
\begin{align*}
\int_V \kappa d \mu & =\int_V  \left((-\Delta)^{s} u\right)  e^{-u}d \mu \\
& = \int_V   \nabla^s  u\nabla^s e^{-u}  d \mu \\
& =\frac{1}{2}\sum_{x \in V}\sum_{y \in  V, \, y \neq x}W_s(x, y)  (u(x)-u(y)) \left(e^{-u(x)}-e^{-u(y)}\right).
\end{align*}
Since $u$ is not a constant,   there is $(u(x_0)-u(y_0))(e^{-u(x_0)}-e^{-u(y_0)}) <0$ for some $x_0\neq  y_0 \in V$,  and then  $(u(x)-u(y))(e^{-u(x)}-e^{-u(y)}) \leq0$ for all $x\neq  y \in V$.
This together with
\eqref{H2}  leads to $\int_V \kappa d \mu <0$, which gives the desired result.

{ \textit{Necessity.}}
 Suppose that $\kappa$ changes sign and $\int_V \kappa d \mu<0.$
Define a function space
\begin{align}\label{lambda-4}
\Lambda_2=\left\{u \in  {W^{s,2}(V)}:  \int_V u d \mu=0 \text{ and }\int_V \kappa e^u d \mu=0\right\} .
\end{align}
 It follows from (\cite{A-Y-Y-2}, Section 4) that  $\Lambda_2 $ is not empty.
We shall solve (\ref{KW-1}) by minimizing the functional
\begin{align}\label{J2}J_2(u)=\frac{1}{2}\int_V| {\nabla^s u}|^{2} d \mu, \ \forall\ { {u \in} \Lambda_2}.\end{align}
Obviously, $J_2$ has a lower bound on   $\Lambda_2$.  Hence this permits us to consider   $  \theta_2=\inf _{ {u \in} \Lambda_2}  J_2(u) $.
Repeating the arguments of the proof of Theorem \ref{T2} ($i$), we can find a function $ {u_0} \in \Lambda_2  $ such that  it is a minimize of $ J_2  (u)$  on $\Lambda_2  $.

The Euler-Lagrange equation of $ {u_0}$ can be calculated   as
 \begin{align*}
\left\{\begin{aligned}
& (-\Delta)^{s}  {u_0}={\lambda_0}  \kappa e^{ {u_0}}, \\
&\lambda_0=\frac{\theta_2}{ \int_V \kappa u_0e^{u_0}d\mu}.
\end{aligned}\right.
 \end{align*}
  We claim that $\lambda_0 \neq 0$. Otherwise, we conclude $ {u_0} \equiv 0 \notin \Lambda_2$ from $(-\Delta)^{s}  {u_0}=0$ and $\int_V  {u_0} d \mu=0$, which is a contradiction. Furthermore,  $\lambda_0>0$   is obtained from  $\int_V \kappa d \mu<0$ and
  \begin{align*}
 \lambda_0 \int_V \kappa d \mu=\int_V e^{- {u_0}} (-\Delta)^{s}   {u_0} d \mu  = \int_V  \nabla^s  u_0 \nabla^s e^{-u_0}d \mu <0.
\end{align*}
And it is not difficult to verify that $u= {u_0}+\ln{ \lambda_0}$ is a desired solution of    (\ref{KW-1}). $\hfill\Box$\\

\textit{\textbf{Proof of Theorem \ref{T2} ($iii$).}}
If $u_0$ is a solution of (\ref{KW}) for some $c_0<0$, namely $e^{-u_0}(-\Delta)^{s} u_0+c_0e^{-u_0}=\kappa   $, then there holds
\begin{align*}
\int_V \kappa d \mu & =\int_V e^{-u_0} (-\Delta)^{s}  u_0 d \mu+c_0 \int_V e^{-u_0} d \mu \\
& = \int_V  \nabla^s  u_0 \nabla^s e^{-u_0}d \mu+c_0 \int_V e^{-u_0} d \mu.
\end{align*}
Similar to the arguments for the sufficiency of Theorem \ref{T2} ($ii$), we obtain $  \int_V  \nabla^s  u_0 \nabla^s e^{-u_0}d \mu <0$. Because of $c_0<0$, we get $c_0 \int_V e^{-u_0} d \mu<0$,
  which gives  the desired result $\int_V \kappa d \mu <0$.$\hfill\Box$

\subsection{Proof of Theorem \ref{T4}}\label{sub4.4}

In this subsection,  we first prove Theorem \ref{T4} ($i$)  by using the method of upper and lower solutions.
If for all $x \in V$,   two functions $u^+$ and $u^-$  satisfy
 \begin{align*}
 (-\Delta)^{s}  u^+(x)-\kappa e^{ u^+(x)}+c \geq 0,\\
 (-\Delta)^{s}  u^-(x)-\kappa e^{ u^-(x)}+c \leq 0,
\end{align*}
then we call $u^+$ and $u^-$  as the upper solution and  lower solution of (\ref{KW}) respectively.
Define an operator $L_{s,\phi}$ acting on a function $u\in C(V)$ as
\begin{align}\label{Lsp}
L_{s,\phi} (u)= \left( (-\Delta)^{s}  +\phi \right)u,\end{align}
where   the function  $\phi  (x)>0$ for all $x \in V$.  Then we have the following

\begin{lemma}\label{L4.1}
Then $L_{s,\phi}$ is bijective. Moreover,   if $ L_{s,\phi}(u)\leq  L_{s,\phi}(v)$, then  $u \leq v$.
\end{lemma}

\begin{proof}

It follows from the H\"{o}lder inequality and the Young  inequality that
\begin{align}\label{3}
\int_{ V} v(-\Delta)^{s} u  d\mu=&  \int_{ V} \nabla^s u \nabla^s v d\mu \notag\\
\leq  &\left(\int_{ V}|\nabla^su|^{2}  d\mu\right)^{\frac{1}{2}}\left(\int_{ V} |\nabla^s v|^{2} d\mu\right)^{\frac{1}{2}}\notag\\
\leq  &{\frac{ 1}{2}}\int_{ V}|\nabla^su|^{2}  d\mu+{\frac{1}{2}}\int_{ V} |\nabla^s v|^{2} d\mu.
\end{align}

 We first prove $L_{s,\phi}$ is injective.
Suppose  $L_{s,\phi}(u)=L_{s,\phi}(v)$, then  there has  $
\phi(v -u ) = (-\Delta)^{s} (u - v). $
This together with \eqref{part} and \eqref{3}   leads to
\begin{align*}
\int_V \phi(v -u )^2d\mu= \int_V (v -u )\left( (-\Delta)^{s} (u - v) \right)d\mu
= 0,
\end{align*}
which implies $u=v$ from $V$ is finite and $\phi>0$.
On the other hand, if we want to  prove $L_{s,\phi}$ is surjective,
 we just need to show that for any given $f \in C(V)$, the equation
  \begin{align}\label{5}(-\Delta)^{s} u +\phi u=f\end{align}
  has a solution $u$. Denote an energy functional $J_3(u)$ as
  \begin{align*}
J_3(u)=\frac{1}{2} \int_V|\nabla^s u|^{2} d \mu+\frac{1}{2} \int_V \phi u^2 d \mu-\int_V f u d \mu, \quad \forall u \in C(V) .
\end{align*}
 Based on the method of Lagrange multipliers, one can calculate the Euler-Lagrange equation of $J_3(u)$ as \eqref{5}.
 By the H\"{o}lder inequality, we get
\begin{align*}
J_3(u) & \geq \frac{1}{2} \int_V \phi u^2 d \mu-\int_V f u d \mu\\
& \geq \frac{1}{2} \min_{x\in V}\phi(x) \| u\|_2^2-\|f\|_2 \| u \|_2  \rightarrow+\infty
\end{align*}
 as $\|u\|_2 \rightarrow+\infty$. Since $C(V)$ is a finite-dimensional linear space,  then $J_3(u)$ attains its minimum by some $u_0\in C(V)$,  which is a solution of   \eqref{5}, and then $L_{s,\phi}$ is bijective.

Moreover, if $ L_{s,\phi}(u)\leq  L_{s,\phi}(v)$, then $ u\leq  v$. Suppose not, then there exists some $x_0\in V$ such that
  \begin{align*} \max _{x \in V}\left(u(x)-v(x)\right)=u \left(x_0\right)-v\left(x_0\right)>0.\end{align*}
And then
\begin{align*}   v\left(x_0\right)-u(x_0)\leq  v(y) -u(y)  , \  \forall y\in V .\end{align*}
Since $ L_{s,\phi}(u)\leq  L_{s,\phi}(v)$, we have
 $ \phi(u-    v)\leq  (-\Delta)^{s} (v-  u) $, and then
  \begin{align*}
  0<  \phi(x_0)(u(x_0)-    v(x_0))
  \leq & (-\Delta)^{s} ( v-u)(x_0)   \\
  =&\frac{1}{\mu(x_0)}\sum_{y \in  V, \, y \neq (x_0)}W_s(x_0, y)  \left( v(x_0)-u(x_0)-(v(y)-u(y)) \right)\\
  \leq & 0.
  \end{align*}
This is a contradiction, and thus   $u \leq v$.
\end{proof}

\begin{lemma}\label{L2}
Let $c<0$. If there have an upper solution $ u^+$ and a lower solution $ u^-$  of (\ref{KW})   with $  u^+\geq u^-$, then (\ref{KW}) has a solution $u$  satisfying $ u^- \leq u \leq  u^+$.
\end{lemma}

\begin{proof}
This is a discrete version of the argument of Kazdan-Warner (\cite{K-W-1}, Lemma 9.3), and  the method of proof carries over to the setting of graphs.

 Set $\kappa _1(x)=$ $\max \{1,-\kappa(x)\}$, so that $\kappa _1(x) \geq 1$ and $\kappa _1 (x)\geq-\kappa(x)$. Let $L_{s,\phi}$ be defined by  \eqref{Lsp} with $\phi(x)=\kappa_1(x) e^{ u^+(x)}>0$.
 We define inductively $u_{i+1}$ as
\begin{align}\label{17}
L_{s,\phi} (u_{i+1})=\phi u_i+\kappa  e^{u_i} -c,
\end{align}
where $u_0= u^+$ and $i=0$, $1$, $2$, $\cdots$.

We claim that
\begin{align}\label{18}
 u^- \leq u_{i+1} \leq u_i \leq \cdots \leq u_1 \leq   u^+.
\end{align}
To prove the above claim, we just need to prove
\begin{align}\label{18-L}
 L_{s,\phi}(u^- )\leq L_{s,\phi} (u_{i+1}) \leq L_{s,\phi}( u_i) \leq \cdots \leq L_{s,\phi}( u_1) \leq L_{s,\phi} ( u^+)
\end{align}
by Lemma \ref{L4.1}.
We will prove this by the method of induction.
It is easy to see
  \begin{align*}
L_{s,\phi} (u_1)-L_{s,\phi}( u^+) = \kappa(x) e^{u^+}-c- (-\Delta)^{s}  u^+ \leq 0 .
\end{align*}
 Suppose $ L_{s,\phi}(u_i)-L_{s,\phi}(u_{i-1})$, then $u_i \leq u_{i-1}$, and we calculate by using the mean value theorem
\begin{align*}
L_{s,\phi}(u_{i+1})-L_{s,\phi}(u_i) 
 & =\phi(u_i-u_{i-1} )+\kappa \left(e^{u_i}-e^{u_{i-1}}\right) \\
& = \kappa_1  e^{ u^+ }(u_i-u_{i-1} )+\kappa e^{\xi}(u_i-u_{i-1} ) \\
& \leq
 \kappa_1(e^{ u^+}-e^{\xi})(u_i-u_{i-1} ) \\
& \leq 0,
\end{align*}
where $u_i < \xi < u_{i-1}$. Then by induction, $$L_{s,\phi} (u_{i+1}) \leq L_{s,\phi} (u_i )\leq \cdots \leq L_{s,\phi} (u_1 )\leq L_{s,\phi}  (u^+)$$
 for any $i$.
 Similarly as above,
we also have by induction $ L_{s,\phi}(u^- )\leq  L_{s,\phi}(u_{i+1})$   for all $i$.
Therefore, \eqref{18-L} holds and the claim in (\ref{18}) follows.

Since $V$ is finite and (\ref{18}), it is easy to see that there exists a function $ u^*$ satisfying $ u^- \leq u^* \leq  u^+$, such that   $u_i \rightarrow u^*$ uniformly on $V$. Passing to the limit $i \rightarrow+\infty$ in  (\ref{17}), one concludes that $u^*$ is a solution of (\ref{KW}).
\end{proof}

We will demonstrate that (\ref{KW}) has an infinite number of lower solutions. This simplifies the proof of Theorem \ref{T4} (i), as we only need to find its upper solution.

\begin{lemma}\label{L3}
If $c<0$,
then for any   $f\in C(V)$,  (\ref{KW})  always has a lower solution $ u^-$ with $  u^-\leq f$.
\end{lemma}

\begin{proof}
Let $u_{\ell} \equiv-{\ell}$ for some positive constant ${\ell}$. Since $V$ is finite, we have
  \begin{align*}
(-\Delta)^s u_{\ell}  -k  e^{u_{\ell}}+c= -k  e^{-\ell}+c \rightarrow c \ \text { as }\ \ell \rightarrow+\infty,
\end{align*}
uniformly with respect to $x \in V$. Noting that $c<0$, we can find a sufficiently large $\ell$ such that $u_{\ell}$ is a lower solution of (\ref{KW}), which can be chosen smaller than any given function $f\in C(V)$.
\end{proof}

\textit{\textbf{Proof of Theorem \ref{T4} ($i$).}}
Define
\begin{align}\label{c(k)}
c_{s,\kappa}=\inf \left\{c<0: (-\Delta)^{s} u=\kappa e^{u}-c  \, \text { has an upper solution in } V \right\} .
\end{align}
Clearly, if $ u^+$ is an upper solution for a given $c<0$, then $ u^+$ is also an upper solution for all $\tilde{c}$ with $c \leq \tilde{c}<0$, which leads to   $ c_{s,\kappa}<0$.
For any $c$ with $c_{s,\kappa}<c<0$, it is easy to see that (\ref{KW}) has at least an upper solution $u^+$. By Lemma \ref{L3}, we can find a lower solution $u^-$ to (\ref{KW}) with $u^- \leq u^+$. Then by Lemma \ref{L2}, (\ref{KW}) has a solution for this $c$.
In addition,  for any $c$ with $c<c_{s,\kappa}$, (\ref{KW}) has no upper solutions by the definition of $c_{s,\kappa}$.
Since every solution is also an upper solution, (\ref{KW}) has no solution for this $c$. In conclusion,
the constant $c_{s,\kappa}$ defined by  \eqref{c(k)} is what we desire.
$\hfill\Box$\\

Through the definition of $c_{s,\kappa}$ in \eqref{c(k)},   $c_{s,\kappa}= -\infty$ means that
(\ref{KW}) is solvable for any $c<0$, and we denote its solution as $u_c$. In the following lemma, we give properties of $u_c$.

\begin{lemma}\label{L1}
Let $c<0$ and $\kappa$ are given by  (\ref{KW}).
Then the equation
\begin{align}\label{8-2}
\left(  (-\Delta )^{s} -c\right) \varphi=-\kappa
\end{align}
has a unique solution  $\varphi_0$.
Moreover, there holds  for any $x\in V$,
\begin{align}\label{9-2}\varphi_0(x) \geq e^{-u_c(x)}.\end{align}
\end{lemma}

\begin{proof}
 Uniqueness of the solution to \eqref{8-2} can be proved following   Lemma \ref{L4.1}  with $\phi=-c>0$.
Also from Lemma \ref{L4.1}, if we want  to prove \eqref{9-2}, we just need to   show
\begin{align}\label{10-2}
\left(  (-\Delta )^{s} -c\right)  \varphi_0(x)\geq  \left(  (-\Delta )^{s} -c\right)   e^{-u_c(x)}.
\end{align}
Using the fact of $1+x  \leq   e^x $  for any $x\in \mathbb{R}$, we derive
$
 e^{-u_c(x)} -e^{-u_c(y)} 
\leq - e^{-u_c(x)}  \left(u_c(x) -u_c(y)\right) .
$
Then there holds
\begin{align*}
 (-\Delta)^se^{-u_c(x)}
=& \frac{1}{\mu(x)}\sum_{y \in  V, \, y \neq x} W_s(x, y) \left(e^{-u_c(x)} -e^{-u_c(y)}\right)\\
\leq & - e^{-u_c(x)} \frac{1}{\mu(x)} \sum_{y \in  V, \, y \neq x} W_s(x, y) \left(u_c(x) -u_c(y)\right)\\
=&- e^{- u_c(x)}  (-\Delta)^s  u_c(x).
\end{align*}
Since $u_c$ is a solution of (\ref{KW}),  then
 \begin{align*}
 \left( (-\Delta )^{s} -c\right)e^{-u_c(x)}
\leq        -e^{- u_c(x)}  (-\Delta)^s  u_c(x)-c e^{- u_c(x)}
 =   -\kappa(x),
\end{align*}
which together with \eqref{8-2} leads to \eqref{10-2}. This gives the desired result.
\end{proof}

\begin{lemma}\label{L4.2}
 For any  $f \in C(V)$, there
exists a unique solution $u$  (up to a constant) to the equation
 \begin{align}\label{7}
(-\Delta)^su = f-\bar{f }.\end{align}
 \end{lemma}

 \begin{proof}
Using a process similar to the proof of Lemma \ref{L4.1}, we can establish the uniqueness of equation \eqref{7}  (up to a constant).
 The details are somewhat tedious and are omitted here.
We now focus on the solvability of   \eqref{7} for any given $f \in C(V)$.  If $f$ equals a constant, then a solution can be chosen to be an arbitrary constant since the right-hand side satisfies $f-\bar{f}=0$ in this case.
 Hence we suppose that $f$  is not a constant in the following.

  Define a function space as $\Lambda_4=\{u\in C(V): \int_V u d \mu=0 \}$.  Obviously,   $\Lambda_4 \neq \varnothing $ from $f-\bar{f} \in \Lambda_4$ for any $f\in C(V)$. Define    an energy functional as
  \begin{align*}
J_4(u)=\frac{1}{2} \int_V|\nabla^s u|^{2} d \mu-\int_V(f-\bar{f})u d \mu, \quad \forall u \in \Lambda_4,
\end{align*}
 Next we estimate $J_4(u)$.
 By the H\"{o}lder inequality and Lemma  \ref{Lpoincare}, we obtain
\begin{align*}
\left|\int_V(f-\bar{f})(u-\bar{u}) d \mu\right| \leq\|f-\bar{f}\|_2\, \|u-\bar{u}\|_{2}
 \leq  C \|\nabla^s u\|_{2},
\end{align*}
where the constant $C>0$ depends only on   $f$ and $G$.
Then there has
\begin{align*}
J_4(u)
=\frac{1}{2} \int_V|\nabla^s u|^{2} d \mu-\int_V(f-\bar{f})(u-\bar{u})d \mu
 \geq \frac{1}{2} \|\nabla^s u\|_{2}^{2} -C \|\nabla^s u\|_{2}.
\end{align*}
From the fact that     if  $x> 4 a  $, then  $ax< x^{2}/4$ for any $a>0$,  we derive that $J_4(u) \geq      \|\nabla^s u\|_{2}^{2}/4$  if  $ \|\nabla^s u\|_{2}$ is big enough.
 Using  Lemma  \ref{Lpoincare},  we know if $u \in \Lambda_4$, then
$
J_4(u)  \rightarrow+\infty
$
 as $\|u\|_{2}\rightarrow+\infty$.
 This implies that $J_4(u) $ attains its minimum at some $u_0 \in\Lambda_4$.
   Moreover, it is not difficult to verify that the Euler-Lagrange equation of $J_4(u)$ is \eqref{7},  and then $u_0 $  is a solution of   \eqref{7}.
  \end{proof}

\textit{\textbf{Proof of Theorem \ref{T4} ($ii$).}}
 { \textit{Sufficiency.}}
 Denote   $u_c$  as a solution of (\ref{KW}) for any fixed $c<0$. It is clear that $\kappa \not \equiv 0$ from $\int_{V} \kappa d \mu<0$. We now prove $\kappa(x) \leq 0$ for all $x \in V$ by contradiction.  Suppose  there exists some $x_0 \in V$, such that $\kappa(x_0)=\max_{x\in V}\kappa(x)>0$. Let  $\varphi_0$  be a  unique
 solution to  \eqref{8-2}.
 One can get
  \begin{align*}\lim_{c\rightarrow -\infty}c\varphi_0(x_0)=\kappa(x_0)>0\end{align*}
 easily by the orthogonally diagonalization method.
This implies $\varphi_0(x_0) < 0$. However, it follows from Lemma \ref{L1} that $\varphi_0 (x_0)\geq e^{-u_c(x_0)} > 0$, which is a contradiction.

 { \textit{Necessity.}}
 We  first let $v$ be a solution of $(-\Delta)^{s} v=\kappa -\bar{\kappa}$.
The existence of $v$ can be seen in  Lemma \ref{L4.2}.
Since $\bar{\kappa}<0$, for any $c<0$,  we can choose a constant $a> \bar{\kappa}/c >0$ and a sufficient constant  $b>  \ln a-a v$. Denote $u^+=a v+b$, then $e^{u^+}>a  $. It follows from $\kappa(x) \leq 0$ for all $x \in V$ that
\begin{align*}
(-\Delta)^{s} u^+- \kappa e^{ u^+} +c& =a  (-\Delta)^{s} v- \kappa e^{u^+}+c \\
& =a  \kappa- a  \bar{\kappa}- \kappa e^{u^+}+c \\
& > \kappa \left(a   -e^{u^+}\right) \\
& \geq 0 .
\end{align*}
Hence $ u^+$ is an upper solution of (\ref{KW}).  From Lemma \ref{L3}, we can find a lower solution $u^-$ to (\ref{KW}) with $u^- \leq u^+$. Then it follows from Lemma \ref{L2} that (\ref{KW}) has a solution.  Consequently, $c_{s,\kappa}=-\infty$, which gives the desired result.
$\hfill\Box$\\

Suppose  $c_{s,\kappa}\neq -\infty$. We can choose a sequence of decreasing numbers $c_i$ such that
$-\infty< c_{s,\kappa}< c_i <c_{s,\kappa}/2<0$    and $ c_i \rightarrow c_{s,\kappa}$ as $ i \rightarrow +\infty .$
The main difficulty for the final proof  is to find a uniformly bounded solution $u_i$ to this equation
\begin{equation}\label{KW-i}
(-\Delta)^{s} u=\kappa e^{u}-c_i \ \text { in } \ V.
\end{equation}
We shall prove the existence of such $u_i$ by the following two lemmas.

\begin{lemma}\label{L4.3}
For any  $c_i$, \eqref{KW-i} has  a solution $u_i$  that satisfies $u_i\geq M$ for some constant $M$.
\end{lemma}
\begin{proof}
On the one hand, for any fixed $c_i$, there is a $\tilde{c}_i$ such that $  c_{s,\kappa}<\tilde{c}_i< c_i$. It follows from the definition of $c_{s,\kappa}$ in  \eqref{c(k)} that there exists some function $u_i^+\in C(V)$ satisfying $(-\Delta)^{s} u_i^+=\kappa e^{u_i^+}- \tilde{c}_i $. Hence we obtain
\begin{align}\label{13}
 (-\Delta)^s  u_i^+-\kappa e^{u_i^+}+c_i = c_i-\tilde{c}_i > 0,
\end{align}
and thus $u_i^+$ is a super solution to   (\ref{KW-i}). On the other hand, since $\int_{V} \kappa d \mu<0$, we can take some $x_0\in V$ such that  $\kappa(x_0)=\min_{x\in V}\kappa(x)<0$. Then it follows from $c_i <c_{s,\kappa}/2$ that
for any  constant  $M$ satisfying $M<\ln  (c_{s,\kappa}/2 \kappa(x_0) )$,  there holds
\begin{align}\label{11}
  (-\Delta)^{s} M-\kappa(x) e^{M}+c_i\leq-\kappa(x_0) e^{M}+c_i<0 .
\end{align}
 Then the constant function $M$ is a lower solution to (\ref{KW-i}).

Claim that $M<u_i^+(x)$.  We just need to prove $M<u_i^+(x_1)=\min _{x \in V} u_i^+(x)$. There is
  \begin{align*}
(-\Delta)^{s} u_i^+(x_1)  =\frac{1}{\mu(x_1)}\sum_{y \in  V, \, y \neq x_1} {W_s(x_1, y)} \left(u_i^+(x_1)-u_i^+(y)\right)\leq 0.
\end{align*}
In view of \eqref{13}, we obtain
  \begin{align*}
\kappa(x_1) e^{u_i^+(x_1)}< (-\Delta)^{s}  u_i^+(x_1)+{c}_i \leq  {c}_i<0,
\end{align*}
and thus $\kappa(x_1)<0$.
Hence it follows from \eqref{13} and \eqref{11} that
$\kappa(x_1) e^M>
c_i>  
\kappa(x_1) e^{u_i^+(x_1)},$
which leads to $\kappa(x_1)  (e^M-e^{u_i^+(x_1)} )> 0.$
 This together with $\kappa(x_1)<0$ implies  $M<u_i^+(x_1)$, and thus our claim $M<u_i^+$ is follows.
Therefore, it follows from Lemma \ref{L2} that (\ref{KW-i}) has at least one solution $u_i$  satisfying $M\leq u_i\leq u_i^+$.
\end{proof}

\begin{lemma}\label{L4.4}
 If $u_i$ is a solution of $\eqref{KW-i}$, then $u_i$  is uniformly bounded.
\end{lemma}

\begin{proof}
From Lemma \ref{L4.3}, we only need to prove  $u_i$ has a uniform upper bound.
 Since $G$ is a finite graph,   for any  $x \in V\setminus \{x \in V: \kappa(x)=0 \}$, there exists  a constant $K>0$ such that
\begin{align}\label{1.6}
\frac{1}{K} \leq |\kappa(x)|  \leq K,\end{align}
and there exists  $x_i \in V$ such that $u_i (x_i)=\max _{x\in V} u_i(x) $. Hence it is easy to know that
$(-\Delta)^{s} u_i (x_i)
\geq 0$, which leads to
\begin{align}\label{2.3}
\kappa \left(x_i\right) e^{u_i (x_i)}-c_i  =(-\Delta)^{s} u_i (x_i)
 \geq 0.
\end{align}
Our proof is closely related to the sign of $\kappa \left(x_i\right)$.   Next, we will discuss three cases that may occur as follows:

We first consider the simplest case of $\kappa (x_i)<0$.  It follows from \eqref{1.6} and \eqref{2.3} that
  \begin{align*}
\frac{1}{K} e^{u_i (x_i)}+c_{s,\kappa} \leq&-\kappa \left(x_i\right) e^{u_i (x_i)}+c_i\\ =&-(-\Delta)^{s} u_i (x_i)\\\leq& 0,
\end{align*}
which leads to
$
u_i(x)
\leq \ln \left(-Kc_{s,\kappa}   \right),
$
and then the proof is done.

We next consider the second case of $ \kappa (x_i)>0$.  It follows from Lemma \ref{L4.3} that $ u_i>M$.
This together with $ c_i<0$,   \eqref{1.6} and \eqref{2.3} leads to
\begin{align*}
\frac{1}{K} e^{u_i (x_i)}&\leq  \kappa \left(x_i\right) e^{u_i (x_i)} \\ &=(-\Delta)^{s} u_i (x_i)+c_i \\
&\leq \frac{1}{\mu(x) } \sum_{y \in  V, \, y \neq x_i}  {W_s(x_i, y)}  \left(u_i^+(x_i)-u_i^+(y)\right)   \\
& \leq C_1(V)\left(u_i(x_i)-M\right),
\end{align*}
where the constant
\begin{align}\label{C1}
C_1(V)=\frac{n}{ \min_{x\in V} \mu(x)}  \max_{y \in  V, \, y \neq x} {W_s(x, y)}  >0 .
\end{align}
Suppose $u_i(x_i)\rightarrow +\infty$ as $ i \rightarrow +\infty $, then there holds
  \begin{align*}  \frac{1}{K} \leq \frac{C_1(V)\left(u_i(x_i)-M \right)}{ e^{u_i (x_i)} }\rightarrow 0 \end{align*}
as $ i \rightarrow +\infty$, which is impossible. Then $u_i$ has a uniform upper bound.

The final case is $ \kappa (x_i)=0$.
Since $(-\Delta)^{s} u_i (x_i)=-c_i >0$, then  $u_i$ is not a constant function and
we have
\begin{align*}
 -c_{s,\kappa}&>-c_i \\
 & =\frac{1}{\mu(x_i)}\sum_{y \in  V, \, y  \neq  x_i} W_s(x_i, y)   \left(u_i^+(x_i)-u_i^+(y)\right)  \\
& \geq   C_2(V) \Big(u_i(x_i)-\min _{y \in  V, \, y  \neq  x_i}u_i(y)\Big),
\end{align*}
where the constant $$C_2(V)= \frac{1} {\max _{x \in  V}\mu(x)}\max _{y \in  V, \, y  \neq  x } {W_s(x , y)}  >0.$$
Hence we conclude that
\begin{align}\label{2.6}
0 \leq u_i (x_i) \leq \min _{y \in  V, \, y  \neq  x_i}u_i(y)-\frac{c_{s,\kappa}}{ C_2(V)}. 
\end{align}
Since  $u_i$ is not a constant function,
there exists  $\tilde{x}_i \in V$ with $ \tilde{x}_i \neq  x_i$ such that $u_i (\tilde{x}_i  )=\min _{x\in V} u_i(x)$.
In the next, we will prove $u_i(\tilde{x}_i  )$ has a uniform upper bound. We may assume that $u_i(\tilde{x}_i)>0$, for otherwise, $u_i(\tilde{x}_i) \leq 0$ as we desired.
Noting that $c_i<0$, we obtain
\begin{align*}
\kappa(\tilde{x}_i) e^{u_i(\tilde{x}_i)}<\kappa(\tilde{x}_i) e^{u_i(\tilde{x}_i)}-c_i &=(-\Delta)^{s} u_i(\tilde{x}_i) \leq 0,
\end{align*}
which implies $\kappa(\tilde{x}_i)<0$. It follows from  $u_i(\tilde{x}_i)>0$, $u_i (x_i)=\max _{x\in V} u_i(x) $ and $ \tilde{x}_i \neq  x_i$ that
\begin{align*}
 -(-\Delta)^{s} u_i(\tilde{x}_i)
&=\frac{1}{\mu(\tilde{x}_i)}\sum_{y \in  V, \, y  \neq  \tilde{x}_i}  {W_s(\tilde{x}_i, y)} \left(u_i(y)-u_i(\tilde{x}_i)\right) \\
& \leq \frac{1}{\mu(\tilde{x}_i)}\sum_{y \in  V, \, y  \neq  \tilde{x}_i}  {W_s(\tilde{x}_i, y)}   u_i(y)  \\
& \leq C_1(V)  u_i (x_i ) ,
\end{align*}
where the constant  $C_1(V)$ is given  by \eqref{C1}.
This together with $\kappa(\tilde{x}_i)<0$, \eqref{1.6} and (\ref{2.6}) leads to
\begin{align*}
 \frac{1}{K}  e^{ u_i(\tilde{x}_i)} &\leq -\kappa(\tilde{x}_i) e^{u_i(\tilde{x}_i)}\\
 & =-(-\Delta)^{s} u_i(\tilde{x}_i)-c_i\\
& \leq   C_1(V) \left( u_i(\tilde{x}_i)-\frac{c_{s,\kappa}}{ C_2(V)}\right)-c_{s,\kappa}.
\end{align*}
Suppose $u_i(\tilde{x}_i)\rightarrow +\infty$ as $ i \rightarrow +\infty $, then there holds
  \begin{align*}
  \frac{1}{K} \leq \frac{C_1(V) \left( u_i(\tilde{x}_i)-\frac{c_{s,\kappa}}{ C_2(V)}\right)-c_{s,\kappa}}{e^{ u_i(\tilde{x}_i)}}\rightarrow 0
  \end{align*}
as $ i \rightarrow +\infty$, which is impossible.   This together with (\ref{2.6}) leads to $u(x_i)\leq C$.
\end{proof}

\textit{\textbf{Proof of Theorem \ref{T4} ($iii$).}}
If $c_{s,\kappa}\neq -\infty$, then it is clear that  \eqref{KW-i}  has a uniformly bounded solution $u_i$ from  Lemma  \ref{L4.3} and Lemma \ref{L4.4}.
Therefore, we can choose a subsequence of $u_i$ such that $u_i$ tends to some $u^* \in C(V)$ as $i \rightarrow+ \infty$. Hence $u^*$ is a solution to (\ref{KW}) when $ c =c_{s,\kappa}$.
$\hfill\Box$

\subsection{Proof of Theorem \ref{T2-m} }\label{S8}

In this section, we focus on the solvability of  (\ref{KW}) when $s>1$.  As for the case of $c > 0$, we can obtain Theorem \ref{T2-m} ($i$) by repeating the arguments of the proof of Theorem \ref{T2} ($i$).  We omit the details and leave the proof to interested readers. If $c=0$, the proofs of Theorems \ref{T2-m} ($ii$) and \ref{T2} ($ii$) are not identical, and we only give the differences as follows:\\

\textit{\textbf{Proof of Theorem \ref{T2-m} ($ii$).}}
Let  the space  $\Lambda_2 $ and the functional  $ J_2  (u)$  be given by  \eqref{lambda-4}  and      \eqref{J2} respectively.
We shall solve (\ref{KW}) with $c=0$ by minimizing $ J_2  (u)$.
 Analogous to the proof of Theorem \ref{T2} ($ii$),  we can find a function $ {u_0} \in \Lambda_2  $ such that  it is a minimize of $ J_2  (u)$  on $\Lambda_2  $. The Euler-Lagrange equation of $ {u_0}$ can be calculated as
 \begin{align*}
\left\{\begin{aligned}
& (-\Delta)^{s}  {u_0}(x)={\lambda_0}  \kappa (x)e^{u_0(x)}, \\
&\lambda_0=\frac{\theta_2 }{\int_V \kappa u_0e^{u_0}d\mu}, \   \theta_2  =\inf _{ {u \in} \Lambda_2  }  J_2  (u) .
\end{aligned}\right.
 \end{align*}
We claim that $\lambda_0 \neq 0$. Otherwise, we conclude $ {u_0} \equiv 0 $ from $(-\Delta)^{s}  {u_0}=0$ and $\int_V  {u_0} d \mu=0$,  then $\int_V \kappa d \mu=0$, which is a contradiction.  Furthermore, we claim that $\lambda_0 >0$. Otherwise, it follows from $\theta_2 > 0$ that $\int_V \kappa u_0e^{u_0}d\mu< 0$. Since $ {u_0} \not \equiv 0 $, $\int_V u_0 d \mu=0$ and $V$ is finite, we can get a vertex  $x_0\in V$ such that $u_0(x_0)=\min_{x\in V}u_0(x)<0$. Then
$$0>\int_V \kappa u_0e^{u_0}d\mu\geq  u_0(x_0)e^{u_0(x_0)} \int_V \kappa d\mu ,$$
which is impossible from $ \int_V \kappa d\mu<0$. Therefore, it is not difficult to verify that $u= {u_0}+\ln{ \lambda_0}$ is a desired solution of  (\ref{KW}) when $c=0$.
$\hfill\Box$

For the case $c < 0$,  the method of upper and lower solutions is not available, and then we use the calculus of variations to prove Theorem  \ref{T2-m} ($iii$).\\

\textit{\textbf{Proof of Theorem  \ref{T2-m} ($iii$).}} Let  the space  $\Lambda_1$ and the functional  $ J_1  (u)$  be given by  \eqref{lambda-1}  and      \eqref{J1} respectively.
Now we minimize the functional $ J_1 (u)$ on $ \Lambda_1 $.  For this purpose,  we next prove that $J_1 $ has a lower bound on $\Lambda_1 $.
Write $v=u-\bar{u}$, then $\bar{v}=0$ and $\|\nabla^{s} v\|_2=\|\nabla^{s} u\|_2$. Hence for any $u\in \Lambda_1 $,  we have
$\int_V  \kappa e^v d \mu= c |V| e^{-\bar{u}}>0 $.
 Since $c<0$ and $\kappa(x)<0$ for all $x \in V$, we have $\max _{x \in V} \kappa(x)<0$. Write $\theta= \max _{x \in V}\kappa(x)/c |V|$, and then
 $ {\kappa}(c |V|)^{-1} \geq\theta >0 .$
Therefore,  there holds
\begin{align}\label{e1-m}
 J_1 (u)\geq \frac{1}{2} \int_V|\nabla^{s} u|^2 d \mu-c   |V| \ln \int_V \theta  d \mu-c   |V| \ln \int_V     e^v d \mu  .
\end{align}
According the Jensen inequality, $ \int_V e^v d \mu \geq e^{\bar{v}} {|V|}={|V|} .$
 Given (\ref{e1-m}),  the above inequality leads to that $J_1$ has a lower bound on $\Lambda_1 $.
Repeating the arguments of the proof of Theorem \ref{T2} ($i$), we can find a function $ {u_0} \in \Lambda_1  $ such that  it is a minimize of $ J_1  (u)$  on $\Lambda_1  $.   It is not difficult to check that (\ref{KW}) is the Euler-Lagrange equation of $ J_1 (u)$. This completes the proof of the theorem.$\hfill\Box$

 \section*{Appendix}
At first, we give the matrix representation of the fractional Laplace operator $(-\Delta)^s$ for any $s\in (0,1)$. Assume that $G = (V, E, \mu, {w})$ is a connected finite graph.
 For convenience, 
  we   rewrite $V=\{x_1, x_2, \cdots, x_n\}$ with $n=\sharp V$.
  Define the  measure matrix  of  $G$ as
          \begin{align*}
  \bm{U}=diag \left(\mu(x_1), \mu(x_2),\cdots, \mu(x_n) \right) .
   \end{align*}
    Set  two matrices
    \begin{align*}
 \bm{A}_s=(A_{s,ij})_{n\times n}  \  \text{ and } \  \bm{D}_s=diag \left(d_s(x_1), d_s(x_2),\cdots, d_s(x_n) \right) ,
    \end{align*}
where
     \begin{align*}
     A_{s,ij}=   \left\{\begin{aligned}
    &-W_s(x_i, x_j),  &\text{ if }  i \neq j,\\
    &  0 , &\text{ if }  i= j,
    \end{aligned}\right.
   \end{align*}
   $W_s(\cdot, \cdot)$ is given by \eqref{H}, and $d_s(x_i)=\sum_{j=1}^nA_{s,ij}$. 
Moreover, we define  the  fractional Laplacian  matrix  of  $G$ as
     \begin{align*}   { \bm{L}_s =  \bm{U}^{-1} (\bm{D}_s-\bm{A}_s)}. \end{align*}
   Each function $u:V\rightarrow\mathbb{R}$ can be viewed as a column vector
   $\bm{u}\in\mathbb{R}^n$ in the sense that      \begin{align*}\bm{u}=\left(u(x_1),u(x_2),\cdots,u(x_n)\right)^T.\end{align*}
It is not difficult to verify that
    \begin{align}\label{M}
     \bm{{L_s}} \bm{u}= \left(\begin{array}{lll}
  (-\Delta)^s u(x_1) \\
  (-\Delta)^s  u(x_2) \\
 \ \ \  \quad \   \vdots\\
  (-\Delta)^s u(x_n)
 \end{array}\right).
   \end{align}
   Assume that      \begin{align*} 0=\lambda_1<\lambda_2\leq \cdots\leq \lambda_{n}\end{align*} are the eigenvalues of the fractional Laplace  operator $-\Delta$, and $\phi_1, \phi_2, \cdots, \phi_n$ are the corresponding orthonormal eigenfunctions. It is clear that   $  \phi_1= {1}/{\sqrt{|V|}} $,
 where $|V|=\sum_{x\in V}\mu(x)$ denotes the volume of $V$.
   From Proposition \ref{T},   $\phi_1, \phi_2, \cdots, \phi_n$ remain    orthonormal eigenfunctions corresponding to the eigenvalues  $  \lambda_1^s,\lambda_2^s, \cdots,\lambda_{n}^s$  of the fractional Laplace operator $(-\Delta)^{s}$.
Write
     \begin{align*}\bm{\phi}_i=({\phi}_i(x_1),{\phi}_i(x_2),\cdots,{\phi}_i(x_n))^T\end{align*}
       as a  column vector for any $i=1,\cdots, n$.
From \eqref{M}, we obtain
    \begin{align}\label{phis}
     \bm{{L_s}} \bm{\phi}_i= 
     \lambda_i^s\bm{\phi}_i.
   \end{align}
And then,   $  \lambda_1^s, \lambda_2^s,\cdots, \lambda_{n}^s$ are also
  the eigenvalues of  matrix  $  {\bm{L}_s} $, and $\bm{\phi}_1, \bm{\phi}_2,\cdots, \bm{\phi}_n$ are the corresponding    orthonormal eigenvectors.
    Here the orthonormality concerns the inner product with respect to the measure $\mu$, namely  for any $i, j = 1, 2, \cdots , n,$
    \begin{align}\label{2}
     \langle \bm{\phi}_i, \bm{\phi}_j\rangle= \bm{\phi}_i^T \bm{U} \bm{\phi}_j =
     \sum_{k=1}^n\phi_i(x_k)\phi_j(x_k) \mu(x_k)  =\left\{\begin{aligned}    &1,&{\rm if }\ \ i=j,\\    &0,&{\rm if }\ \ i\not=j.    \end{aligned}\right.    \end{align}
Define  two matrices $ \Phi=  \left(\bm{\phi}_1, \bm{\phi}_2,\cdots, \bm{\phi}_n\right) $ and  $  \Lambda^s=diag (\lambda_1^s, \lambda_2^s,\cdots, \lambda_n^s  ).$   It is not hard to verify that
    \begin{align}\label{Phi}
   \Phi^{T} \bm{U} \Phi=\bm{E}   \ \text{ and }\   
  \Phi^{-1}=\Phi^{T} \bm{U} .
\end{align}
From  \eqref{phis}, we derive
         \begin{align*}
   \bm{L}_s\Phi=( \bm{{L_s}} \bm{\phi}_1,  \bm{{L_s}} \bm{\phi}_2,\cdots, \bm{{L_s}} \bm{\phi}_n)=
    ( \lambda_1^s\bm{\phi}_1,\lambda_2^s\bm{\phi}_2,\cdots,\lambda_n^s\bm{\phi}_n)=  \Phi\Lambda^s,
   \end{align*}
    and    then   the matrix representation of the fractional Laplace operator is
    \begin{align}\label{X} \bm{L}_s=\Phi \Lambda^s \Phi^{-1}  . \end{align}

 Secondly, we prove Proposition \ref{sT1}     through standard techniques of matrix theory. \\

 \textit{\textbf{Proof of  Proposition \ref{sT1}.}} On the one hand,  we set  $  \Lambda =diag (\lambda_1, \lambda_2,\cdots, \lambda_n  ).$ It is well known that
  the matrix representation of the Laplace operator $ -\Delta$ is
     \begin{align*}\bm{L} =\Phi \Lambda  \Phi^{-1}  . \end{align*}
      Therefore, it is obvious that
     \begin{align*}\lim_{s\rightarrow 1^-}\bm{L_s} = \bm{L} ,\end{align*}
from \eqref{X}, and then $(i)$ follows.
On the other hand,   we obtain
\begin{align*}
\lim_{s\rightarrow 0^+}\bm{L_s}&=\Phi
\left(\begin{array}{lllll}
0&0&\cdots&0\\
0&1&\cdots&0\\
   \vdots&\vdots&\ddots&\vdots\\
0&0&\cdots&1
   \end{array}\right) \Phi^{-1}
    =
   \Phi\left(\bm{E}-
\left(\begin{array}{lllll}
1&0&\cdots&0\\
0&0&\cdots&0\\
   \vdots&\vdots&\ddots&\vdots\\
0&0&\cdots&0
   \end{array}\right)
   \right)\Phi^{-1}
   =\bm{E}-     \bm{\phi}_1  \bm{\phi}_1^T \bm{U}
    \end{align*}
  from \eqref{Phi}.
In view of $  \phi_1= {1}/{\sqrt{|V|}} $ and  \eqref{2}, we get
  \begin{align*}
  \bm{\phi}_1=
  \frac{1}{\sqrt{|V|}}(1,1,\cdots,1)^T
  \ \text{ and }\ \sum_{k=1}^n \phi_i(x_k) \mu(x_k)
=0,  \ \forall  i\neq 1.
\end{align*}
 Therefore,  if $u\in C(V)$ satisfying $\int_V{u}d\mu=0$,  then  setting
  $\bm{u}=(u(x_1),u(x_2),\cdots,u(x_n))^T$,    we derive that
\begin{align*}
\lim_{s\rightarrow 0^+}\bm{L_s}\bm{u}  =\bm{u}-
     \frac{1}{|V|}
\left(\begin{array}{lllll}
\mu(x_1)&\mu(x_2)&\cdots&\mu(x_n)\\
\mu(x_1)&\mu(x_2)&\cdots&\mu(x_n)\\
 \ \ \  \vdots&\ \ \ \vdots&\ddots&\ \ \ \vdots\\
\mu(x_1)&\mu(x_2)&\cdots&\mu(x_n)
   \end{array}\right) \left(\begin{array}{lll}u(x_1)\\
   u(x_2)\\
   \ \ \ \vdots\\
   u(x_n)\end{array}\right)
=\bm{u}.
 \end{align*}
 This completes the proof.
$\hfill\Box$

\section*{Acknowledgements}
The authors thank Professor Shuang Liu for many helpful discussions on the heat kernel on finite graphs.
The first named author is supported by the National Postdoctoral Fellowship Program (Grant Number: GZC20231343), and
the second named author is supported by the National Natural Science Foundation of China (Grant Number: 12071245).

\section*{Conflict of interest} 	
On behalf of all authors, the corresponding author declares that there are no conflicts of interest regarding the publication of this paper.


\end{document}